\documentclass{article}
\usepackage{amsmath,amssymb,amsthm, amsfonts}
\usepackage{etaremune, enumerate, float, verbatim}
\usepackage{algorithm}
\usepackage{algorithmic}
\usepackage{url}
\usepackage{hyperref}
\usepackage{graphicx}
\usepackage{url}
\usepackage[mathlines]{lineno}
\usepackage{dsfont} 
\usepackage{mathrsfs} 
\usepackage{soul} 
\usepackage{color}

\definecolor{greene}{rgb}{0,0.7,0}

\definecolor{purpl}{rgb}{.6,0,.6}

\usepackage{tikz}

\numberwithin{figure}{section}   
\numberwithin{table}{section}
\numberwithin{equation}{section} 

\newtheorem{thm}{Theorem}[section]
\newtheorem{cor}[thm]{Corollary}
\newtheorem{lem}[thm]{Lemma}
\newtheorem{prop}[thm]{Proposition}

\newtheorem{obs}[thm]{Observation}

\newtheorem{quest}[thm]{Question}

\theoremstyle{definition}
\newtheorem{rem}[thm]{Remark}

\theoremstyle{definition}
\newtheorem{defn}[thm]{Definition}

\theoremstyle{definition}
\newtheorem{ex}[thm]{Example}

\DeclareMathOperator{\cl}{cl}
\newcommand {\cls} {\cl_-}
\newcommand {\du} {\,\sqcup\,}

\newcommand{\B}{\mathcal{B}}

\newcommand{\bit}{\begin{itemize}}
\newcommand{\eit}{\end{itemize}}
\newcommand{\ben}{\begin{enumerate}}
\newcommand{\een}{\end{enumerate}}
\newcommand{\beq}{\begin{equation}}
\newcommand{\eeq}{\end{equation}}
\newcommand{\bea}{\begin{eqnarray*}}
\newcommand{\eea}{\end{eqnarray*}}
\newcommand{\bean}{\begin{eqnarray}}
\newcommand{\eean}{\end{eqnarray}}
\newcommand{\bpf}{\begin{proof}}
\newcommand{\epf}{\end{proof}\ms}
\newcommand{\bmt}{\begin{bmatrix}}
\newcommand{\emt}{\end{bmatrix}}
\newcommand{\ms}{\medskip}
\newcommand{\noi}{\noindent}
\newcommand{\beqa}{\begin{array}}
\newcommand{\eeqa}{\end{array}}
\newcommand{\ol}{\overline}
\newcommand{\lc}{\left\lceil}
\newcommand{\rc}{\right\rceil}
\newcommand{\lf}{\left\lfloor}
\newcommand{\rf}{\right\rfloor}
\newcommand{\lp}{\left(}
\newcommand{\rp}{\right)}

\newcommand{\wh}{\widehat}

\newcommand{\Xr}{X}

\newcommand{\Z}{\operatorname{Z}}
\newcommand{\Zs}{\operatorname{Z}_-}
\newcommand{\Zp}{\operatorname{Z}_+}

\newcommand{\Dctwo}{\gamma^c_2}
\newcommand{\DtTwo}{\gamma^t_2}
 
\newcommand{\vc}{\tau}
\newcommand{\zbar}{\ol\Z}
\newcommand{\zsbar}{\ol\Z_-}
\newcommand{\zpbar}{\ol\Z_+}

\newcommand{\dbar}{\ol\gamma}

\newcommand{\vcbar}{\ol\vc}
\newcommand{\ibar}{\underline{\alpha}}

\newcommand{\xbar}{\ol X}

\newcommand{\xtar}{\mathfrak{X}}

\newcommand{\ztar}{\mathfrak{Z}}
\newcommand{\vctar}{\mathfrak{C}}
\newcommand{\zirtar}{\mathfrak{Zir}}
\newcommand{\zpirtar}{\mathfrak{Z_+ir}}
\newcommand{\vcirtar}{\mathfrak{Cir}}
\newcommand{\xirtar}{\mathfrak{Xir}}

\DeclareMathOperator{\dir}{dir}
\DeclareMathOperator{\DIR}{DIR}
\DeclareMathOperator{\zir}{zir}
\DeclareMathOperator{\ZIR}{ZIR}
\DeclareMathOperator{\zpir}{z_+ir}
\DeclareMathOperator{\ZpIR}{Z_+IR}
\newcommand{\zsir}{{\rm z}_-{\rm ir}}
\newcommand{\ZsIR}{\Zs\!{\rm IR}}

\DeclareMathOperator{\xir}{xir}
\DeclareMathOperator{\XIR}{XIR}
\DeclareMathOperator{\VCIR}{VCIR}
\DeclareMathOperator{\vcir}{vcir}
\newcommand{\ytar}{\mathfrak{Y}}

\definecolor{purp}{rgb}{.5,0,.5}

\title{Universal perspectives on irredundance for $X$-set parameters}
\author{Bryan Curtis\thanks{
Air Force Research Laboratory, Wright-Patterson Air Force Base, OH 45433, USA
}\and
Mary Flagg\thanks{Department of Mathematics, Statistics and Computer Science, University of St. Thomas, Houston, TX 77006, USA (flaggm@stthom.edu).} \thanks{Corresponding Author.}\and
   Leslie
Hogben\thanks{American Institute of Mathematics, Pasadena, CA 91125, USA
(hogben@aimath.org), Department of Mathematics, Iowa State University,
Ames, IA 50011, USA, and Department of Mathematics, Purdue University, West Lafayette, IN 47907, USA.}}

\begin{document}

\maketitle 

\begin{abstract}  Universal definitions of irredundance for $X$-set parameters are presented using blocking sets.  This approach is modeled on (domination) irredundance (which uses closed neighborhoods as blocking sets) and zero forcing irredundance (which uses forts).  {Results include a chain of inequalities between irredundance parameters and original parameters and the isomorphism theorem for TAR reconfiguration graphs of many irredundance parameters. These results are then applied to PSD forcing irredundance and vertex cover irredundance; the chain of inequalties also applies to  skew forcing irredundance. The upper vertex cover irredundance number becomes part of the Domination Chain used in the study of (domination) irredundance.}  Based on the propagation processes involved in forcing, an alternate universal theory of irredundance is developed using closure operators.
\end{abstract}

 \noi {\bf Keywords} irredundance; $X$-set;   PSD zero forcing; skew zero forcing; vertex cover; closure operator  

\noi{\bf AMS subject classification}  05C69, 05C50, 05C57, 68R10

\section{Introduction}

 The study of irredundance in relation to domination was introduced in 1978 by Cockayne, Hedetniemi and Miller \cite{CHM78} as part of a study of minimal dominating sets and has been studied extensively since then (see \cite{MR21} and the references therein). 
Variations of irredundance such as total irredundance in graphs (e.g., \cite{FHHHK02}) and edge in irredundance hypergraphs (e.g., \cite{CKMU19hyper}) have been studied; irredundance has also been studied more more generally, such as irredundant sets of generators  in Boolean algebras \cite{T93boolean}.
Using domination irredundance as a model,  zero forcing irredundance was defined by Curtis, Hogben, and Roux in \cite{ZIR}.  
In this work we introduce two universal perspectives on irredundance that generalize results for domination irredundance and zero forcing irredundance. The first universal approach allows us to provide a single proof of each result (see Section \ref{s:irred}) and then apply results  to the $X$-irredundant parameters  PSD forcing  number, skew forcing  number, and vertex cover number (see Section \ref{ss:block1x}).

In order to explain the first universal perspective, we define domination and zero forcing irredundance because they serve as the models.  A set $S$ is a \emph{dominating  set} of a graph $G$ if every vertex of $G$ is in $S$ or is a neighbor of a vertex in $S$.  The \emph{domination number} of  $G$ is $\gamma(G)=\min\{|S|: S$ is a 
dominating set of $G\}$.  The \emph{upper domination number} of  $G$ is $\dbar(G)=\max\{|S|: S\mbox{ is a minimal dominating set of }G\}$.  Although the domination number is the more important of these two parameters, the upper domination number (introduced in \cite{CHM78}) is also of interest when all  minimal dominating sets are studied, e.g., in reconfiguration \cite{HS14}.
A set $S$ of vertices in a graph $G$ is \emph{domination irredundant}\footnote{We use the term `domination irredundance' and `DIr-set' 
to avoid confusion because we use the term `irredundance' more abstractly, although in the literature `irredundance' and `Ir-set' are used for these terms.} if each vertex in $S$ dominates a vertex of G that is not dominated by any other vertex in $S$  (this may be a neighbor or itself), called a \emph{private neighbor}.  
The  \emph{upper DIr number} is $\DIR(G) = \max\{|S|: S \text{ is a maximal DIr-set}\}$, and the \emph{lower DIr number} is $\dir(G) = \min\{|S| : S \text{ is a maximal DIr-set}\}$.

 Forts, here called standard forts, play a fundamental role in the study of standard zero forcing and were introduced in \cite{FH18}.  Let $G$ be a graph. A \emph{standard fort} of $G$ is a nonempty set $F\subseteq V(G)$ such that $|F\cap N(v)|\ne 1$ for all $v\in V(G)\setminus F$.  Standard forts were used in \cite{ZIR} to define standard zero forcing irredundance. 
 For $S\subseteq V(G)$ and $u\in S$, a fort $F$ of $G$ is a \emph{private standard fort} of $u$ relative to $S$ provided that $S\cap F = \{u\}$. The set $S\subseteq V(G)$ is a \emph{$\Z$-irredundant set} or \emph{ZIr-set} provided every element of $S$ has a private standard fort, the  \emph{upper ZIr number} is $\ZIR(G) = \max\{|S|: S \text{ is a maximal ZIr-set}\}$, and the \emph{lower ZIr number} is $\zir(G) = \min\{|S| : S \text{ is a maximal ZIr-set}\}$.

Note that both domination and zero forcing irredundance use what we call \emph{blocking sets}: A set $S$ is not a dominating  set (respectively, zero forcing set) of $G$ if it doesn't intersect  the closed neighborhood of every  vertex (respectively, every fort).   In Section \ref{s:irred}, we generalizes  this approach by  defining blocking sets  for parameters that are defined as the minimum cardinality among sets having certain properties, called $X$-sets (see Section \ref{ss:uniXY} for more detail).  Conditions on the blocking sets (satisfied by natural blocking sets for many parameters) are used to define an $X$-irredundant blocking family and associated parameters (Definition \ref{d:Xirred}).

  One example of a result that can be established universally  is the natural chain of inequalities between a parameter and its associated irredundance parameter. For domination this is $\dir(G)\le \gamma(G)\le \dbar(G)\le \DIR(G)$, which is generalized in Corollary \ref{c:param-interlace}.  
 It is shown in Sections \ref{ss:PSD-irred} and \ref{ss:skew-ired} that all four of the parameters are distinct for PSD and skew forcing irredundance (as they are for domination and zero forcing irredundance).  However, the lower vertex cover irredundance number is equal to the vertex cover number (Corollary \ref{c:vcir=tau}).   The Domination Chain \cite{CHM78, MR21} includes the independence number and lower independence number in addition to the upper and lower irredundance numbers and domination and upper domination numbers; see   Equation \eqref{eqDomChain} and the discussion there.
 In Remark \ref{r:ExtDomChain}, we extend the Domination Chain to include the upper vertex cover irredundance number.

  The  TAR Graph Isomorphism Theorems from \cite{XTARiso} and \cite{XTAR-YTAR} are extended to the upper $X$-set irredundance parameter $\XIR$ in Theorem \ref{t:xirred-main} when $X$ is a robust $X$-set parameter  and the $X$-blocking sets satisfy additional conditions (Token Addition and Removal reconfiguration (TAR) graphs are defined in Section \ref{ss:TAR} and robust $X$-set parameters are defined in Section \ref{ss:uniXY}).  This theorem applies to  PSD forcing irredundant sets and vertex cover irredundant sets (but not to skew forcing irredundant sets). 
  as discussed in the sections covering those parameters.
  
 In Section \ref{ss:block1x}, numerous additional parameter-specific results are also established  for PSD forcing, skew forcing, and vertex cover irredundance.   For each of these  upper irredundance numbers, there is a bound related to the domination number (Propositions \ref{p:2dom-PSD}, \ref{p:2dom-skew}, and \ref{p:VCIRmax}), although the actual bound varies.  One or more extreme values of the  upper and lower irredundance numbers  for 
 these parameters is characterized,  including the fact that the  lower PSD forcing irredundance number is one if and only if the graph is a tree (Corollary \ref{c:zpir-tree}).  The values of the irredundance parameters are also determined for well-known graphs such as complete graphs, complete bipartite graphs, paths, and cycles.

The second universal approach to irredundance, described in Section \ref{ss:block2}, uses closure operators to form blocking sets and  is perhaps more natural for propagation parameters such as variants of zero forcing.  We show that a closure operator for $X$ with certain basic properties defines an $X$-irredundasnt blocking family as defined in Section \ref{s:irred}.

\subsection{Notation and terminology}\label{ss:term}

  A graph is a pair $G=(V(G),E(G))$ with vertex set $V(G)$ and edge set $E(G)$; all graphs are simple, undirected and finite with a nonempty vertex set. An edge is a two element subset of $V(G)$, and the edge $\{v,w\} \in E(G)$ can be written as $vw$. The subgraph induced by $S \subseteq V(G)$ will be denoted $G[S]$. For simplicity, $G[V(G)\setminus S] = G-S$. If $S = \{v\}$,  we further simplify this to $G-S = G-v$.

Two vertices $v$ and $w$ are called \emph{adjacent or neighbors} if $vw \in E(G)$. The \emph{open neighborhood} of $v \in V(G)$ is $N(v)=\{w:vw \in E(G)\}$. The \emph{closed neighborhood} of $v$ is $N[v]=N(v) \cup \{v\}$ and $N[S]=\cup_{x\in S} N[x]$ for $S\subseteq V(G)$. The \emph{degree} of the vertex $v$ is $\deg(v)=|N(v)|$. A leaf is a vertex of degree $1$. The  minimum degree of vertices in $G$ is denoted by $\delta(G)=\min\{\deg(v):v \in V(G)\}$. Two vertices $v$ and $w$ are called twins if $N(v)=N(w)$ and adjacent twins if $N[v]=N[w]$. 

Given distinct vertices $v_0$ and $v_k$ in $V(G)$,  a \emph{path} of length $k$  from $v_0$ to $v_k$ is a sequence of distinct vertices $(v_0,\ldots, v_k)$ such that 
$v_i$ is a neighbor of $v_{i+1}$ for every integer $i$, $0\leq i\leq k -1$; a cycle $(v_0, \dots, v_k)$ includes the path from $v_0$ to $v_k$ and also the edge $v_kv_0$. A graph is \emph{connected} if there exists a path between any two distinct vertices. Maximal connected induced subgraphs of a disconnected graph are call \emph{connected components}.   The notation $A\du B$ indicates that sets $A$ and $B$ are disjoint, or for   graphs $G$ and $H$, 
the notation $G\du H$ means that the graphs are \emph{disjoint}, i.e., $V(G)\cap V(H)=\emptyset$. 
A \emph{cut-vertex} is a vertex $v \in V(G)$ such that $G-v$ has more connected components than $G$. The \emph{join} of disjoint graphs $G$ and $H$ is the graph $G \vee H$ with vertex set $V(G \vee H)=V(G) \du V(H)$ and edge set $E(G) \du E(H) \du \{ gh:g \in V(G), h \in V(H)\}$.

Given an integer $n \geq 1$, a  complete graph $K_n$ of order $n$ has every possible edge between its vertices. 
Let $\ol{K_n}$ be the graph with $n$ vertices and no edges. 
Given an integer $n \geq 1$, a  path graph $P_n$ (or a cycle graph $C_n$ for $n\ge 3$) is a graph of order $n$ with $V(P_n)=V(C_n)=\{v_1,\dots,v_n\}$ and $E(P_n)=\{v_iv_{i+1}:i=1,
\dots,n-1\}$ (or $E(C_n)=\{v_iv_{i+1}:i=1,
\dots,n-1\}\cup\{v_nv_1\}$).
Given  disjoint sets  of vertices $A=\{a_1,\dots,a_p\}$ and $B=\{b_1,\dots,b_q\}$, the complete bipartite graph   $K_{p,q}$ has $V(K_{p,q})=A\du B$ and $E(K_{p,q})=\{a_ib_j: i=1,\dots,p, j=1,\dots,q\}$; we assume $1\le p\le q$.

\subsection{Universal approach to $X$- and $Y$-sets}\label{ss:uniXY}
A universal approach has been introduced in several contexts related to variants of zero forcing and power domination. For reconfiguration graphs this approach was introduced in \cite{PDrecon} and \cite{XTARiso} and extended in \cite{XTAR-YTAR}. 
We use definitions from \cite{XTAR-YTAR}.  Let $W$ be a property defined on subsets of the vertex set of each graph. If for every graph isomorphism $\varphi: V(G) \to V(G')$, $S\subseteq V(G)$ has property $W$ implies $\varphi(S)$ has property $W$, then $W$ is a \emph{vertex-set property} $W$ and $S\subseteq V(G)$ is a \emph{$W$-set} if it has property $W$. 
 A \emph{cohesive} vertex-set property $W$ requires that every graph has at least one $W$-set, and a well-defined graph parameter for which the value associated to each graph is solely determined by a cohesive property is called a \emph{cohesive parameter}.

A \emph{super $X$-set 
parameter}  is a cohesive parameter $X$ such that $S$ is an $X$-set and $S\subseteq S'$ implies  $S'$ is an $X$-set, and  $X(G)$ is  the minimum cardinality of an $X$-set of $G$. The \emph{upper $X$ number}, denoted by  $\xbar(G)$, is the maximum cardinality of a minimal $X$-set.   Note that $V(G)$ is an $X$-set for any graph $G$ and any super $X$-set parameter (because $G$ has some $X$-set and the $X$ has the superset property).
 A \emph{sub $Y$-set parameter}  is a cohesive parameter $Y$ such that $S$ is a $Y$-set and $S'\subseteq S$ implies  $S'$ is a $Y$-set, and  $Y(G)$ is  the maximum cardinality of a $Y$-set of $G$.   Note that $\emptyset$ is a $Y$-set for any graph $G$ and any sub $Y$-set parameter. 

If one plans to study disconnected graphs, it is useful to require component consistency: A super $X$-set  parameter is \emph{component consistent} if $S$ is an $X$-set  of $G$ if and only if $S\cap V(G_i)$ is an $X$-set of $G_i$ for $i=1,\dots,k$ where   $G_i,i=1,\dots,k$ are the connected components of $G$  and similarly for a $Y$-set parameter.  
The domination number, power domination number,  (standard, PSD, and skew) zero forcing numbers, and the vertex cover number are all examples of component consistent super $X$-set parameters (and all are robust, as defined in the next paragraph). The domination irredundance number and the  standard zero forcing irredunundance number are examples of component consistent  sub $Y$-set parameters (and both are robust, as defined in the next paragraph).

An additional condition is needed to obtain many universal results in \cite{PDrecon, XTARiso,XTAR-YTAR}. A \emph{robust  $\Xr$-set parameter} is a component consistent  super $X$-set parameter $X$  such that   the $X$-sets of $G$ satisfy the following condition:  
\ben[]
\item ($(n-1)$-set)
If $G$ is a connected graph of order  $n\ge 2$, then every set of $n-1$ vertices is an $\Xr$-set.
\een 
A \emph{robust  $Y$-set parameter} is a component consistent  sub $Y$-set parameter $Y$  such that   the $Y$-sets of $G$ satisfy the following condition  (a \emph{singleton} is a set with exactly one element):  
\ben[]
\item 
(Singletons)
If $G$ is a connected graph of order  at least two,
then every singleton is a $Y$-set. 
\een  

The definition of robust $X$-set parameter in \cite{XTAR-YTAR} is adapted from those in \cite{PDrecon, XTARiso}, whereas the definition of robust $Y$-set parameter was introduced in \cite{XTAR-YTAR} motivated by the following connection observed there between $X$- and $Y$-set parameters:  Let $G$ be a graph.
For a  sub $Y$-set parameter $Y$, a subset $S\subseteq V(G)$ is an $X_Y$-set if and only if $V(G)\setminus S$ is a $Y$-set of $G$, and  the parameter $X_Y(G)$ is the minimum cardinality of an $X_Y$-set of $G$.
For a super $X$-set parameter $X$, a subset $S\subseteq V(G)$ is a $Y_X$-set if and only if $V(G)\setminus S$ is an $X$-set of $G$, and  the parameter $Y_X(G)$ is the maximum cardinality of  a $Y_X$-set of $G$.
The parameters $Y$ and $X_Y$ (or $X$ and $Y_X$) are called \emph{base graph complement  parameters}. 
This connection allows translation of $X$-set parameter results to $Y$-set parameters.

\subsection{Parameters}\label{ss:param}

 Standard zero forcing was introduced in multiple applications, including control of quantum systems and as an upper bound to the maximum nullity of a real symmetric matrix $ A=[a_{ij}]$ whose off-diagonal entries $a_{ij}, i\ne j$ are nonzero or zero according as the graph has  edge $ij$.  Other types of zero forcing were defined to serve as similar bounds for maximum nullity among positive semidefinite or skew-symmetric matrices (for skew forcing, the distinguishing feature is that the diagonal entries must all be zero; PSD forcing involves deeper linear algebraic results).

Starting with an initial set of blue vertices and the remaining vertices white, a zero forcing process   applies  a \emph{color change rule} repeatedly to color white vertices blue; the color change rule determines the type of zero forcing.  Here we discuss three types of zero forcing: standard, which was introduced in \cite{AIM08}; positive semidefinite (PSD), which was introduced in \cite{smallparam}; and skew, which was introduced in \cite{IMA10}.    
\bit 
\item {\bf Standard color change rule:} Any blue vertex $u$ can change the color of a white vertex $w$ to blue if $w$ is the only white neighbor of $u$. 
\item {\bf PSD color change rule:}  Let $B$ be the set of blue vertices and let $W_1,\dots, W_k$ be the sets of vertices of the $k\ge 1$ components of $G-B$.  If $u\in B$, $w\in W_i$, and $w$ is the only white neighbor of  $u$ in $G[W_i\cup B]$, then $u$ can change the color  $w$ to blue. 
\item {\bf Skew color change rule:} Any vertex $u$ can change the color of a white vertex $w$ to blue if $w$ is the only white neighbor of $u$. 
\eit
 Given a graph $G$, set $S\subseteq V(G)$, and a  color change rule, the \emph{final coloring} is the set of blue vertices that result from applying the color change rule until no more changes are possible.
A \emph{(standard, PSD, or skew) zero forcing set}   for  $G$  is a subset of vertices $S$ such that  $V(G)$ is the final coloring of $S$.  
The \emph{(standard, PSD, or skew)  zero forcing number} of a graph $G$ is  the minimum cardinality of a standard, PSD, or skew zero forcing set; this parameter is denoted by $\Z(G)$, $\Zp(G)$, or $\Zs(G)$, respectively. All these parameters are also discussed in \cite[Chapter 9]{HLS22book}.  The  \emph{upper zero forcing number} $\zbar(G)$ of a graph $G$ has been defined to be the maximum cardinality of a minimal zero forcing set of $G$ \cite{XTARiso}; the  \emph{upper PSD forcing number} $\zpbar(G)$ and the \emph{upper skew forcing number} $\zsbar(G)$ are defined analogously \cite{XTAR-YTAR}.

A set $S$ of vertices in a graph $G$ is \emph{independent} (or is an \emph{independent set}) of $G$ if no two distinct vertices in $S$ are adjacent.  The \emph{independence number} of $G$, denoted by $\alpha(G)$, is the maximum cardinality of an independent set of vertices of $G$ and the \emph{lower independence number} $\underline{\alpha}(G)$\footnote{In the literature, $\underline{\alpha}(G)$ is often dented by $i(G)$.} is the minimum cardinality of a maximal independent set of $G$. 
A set $S$ of vertices in a graph $G$ is a \emph{vertex cover} of $G$ if every edge of $G$ has at least one of its endpoints in $S$.  
The \emph{vertex cover number} of $G$, denoted by $\vc(G)$, is the minimum cardinality of a vertex cover  of $G$, so the upper vertex cover  number of $G$ is denoted  by $\vcbar(G)$.
It is
well-known (and easy to see) that a set $S\subseteq V(G)$ is independent if and only if $V(G)\setminus S$ is a vertex cover of $G$.  
The independence number and vertex cover number are widely studied parameters and the lower independence number has been studied in conjunction with the Domination Chain.

\section{Blocking sets and irredundance}\label{s:irred}

In this section we develop a universal framework for irredundance parameters defined from super $X$-set parameters   using known definitions and results for domination irreundance and zero forcing irredundance as models.    This universal perspective  provides natural definitions of irredundance for related $X$-set parameters for which irredundance has not yet been studied, including PSD and skew forcing and vertex covering.  
We begin with a very general definition of blocking set.  We then identify additional properties, such as component consistency and $X$-irredundance (see Section \ref{ss:irred}),
that are useful for obtaining results about irredundance parameters; domination irredundance and zero forcing irredundance have these properties. 
  In Section \ref{ss:TAR} we use the properties and additional conditions to establish a TAR graph isomorphism theorem for irredundance parameters analogous to that proved  in \cite{XTAR-YTAR}.

\begin{defn}\label{def:X-block-fam}
Let $X$ be a super $X$-set parameter. For a graph $G$, a subset of vertices $R\subseteq V(G)$ is an  \emph{$X$-blocking set} if $V(G)\setminus R$ is not an $X$-set, and $R$ is a minimal $X$-blocking set if $R$ is an $X$-blocking set and $R'\subsetneq R$ implies $R'$ is not an $X$-blocking set.
An \emph{$X$-blocking family} $\B_X$ is defined by assigning to each graph $G$ a set $B_X(G)$ of (zero or more) $X$-blocking sets of $G$ with the property that if $G$ and $H$ are isomorphic then there exists a relabeling of $V(G)$ such that $B_X(G) = B_X(H)$.  
\end{defn}

 A super $X$-set parameter usually has more than one $X$-blocking family. However, for specific parameters, there is often a natural choice for a blocking family. For example, forts act as natural blocking sets for zero forcing and the closed neighborhood of a vertex acts as a natural blocking set for domination. We fix the notation $B_D(G)=\{ N[v]: v\in V(G)\}$ and $B_{\Z}(G)$ to be the set of all forts of $G$  \cite{ZIR},  thus also fixing $\B_D$ and $\B_{\Z}$.

\begin{defn}
The \emph{maximal $X$-blocking family}, denoted by $\B_X^{\max}$, is formed by assigning the set of all $X$-blocking sets of $G$, denoted $B_X^{\max}(G)$,  to each graph $G$. The \emph{minimal $X$-blocking family}, denoted by $\B_X^{\min}$, is formed by assigning the set of all minimal $X$-blocking sets of $G$, denoted $B_X^{\min}(G)$, to each graph $G$.  The \emph{empty $X$-blocking family}, denoted by $\B_X^{\emptyset}$, is formed by defining $B_X^{\emptyset}(G)=\emptyset$ for each graph $G$.
\end{defn}

  The minimal $\Z$-blocking family has $B_{\Z}^{\min}(G)$ equal to the set of minimal forts of $G$.  The maximal $\Z$-blocking family has $B_{\Z}^{\max}(G)$ equal to the set of supersets of forts. 
 The maximal $D$-blocking family for domination has all supersets of the closed neighborhoods as $B^{\max}_D(G)$ and $\B^{\min}_D$ has the set of minimal closed neighborhoods of $G$ as $B^{\min}_D(G)$ ($N[w]$ is a \emph{minimal closed neighborhood} if $N[v]\subseteq N[w]$ implies $v=w$ for all $v\in V(G)$).
 
 Although it does not happen with domination or standard zero forcing, for some parameters such as skew forcing and vertex covering, there exist graphs $G$ for which the empty set is an $X$-set. 
 In this case there are no $X$-blocking sets for $G$ and necessarily $B_X(G)=\emptyset$.

The next definition provides notation for the slightly different interpretations of maximal and minimal in `maximal $X$-blocking family' and `minimal $X$-blocking family'.

\begin{defn} Suppose $X$ is a super $X$-set parameter and  $\B_X$ and $\B'_X$ are both $X$-blocking families.  
The notation $\B_X\subseteq \B'_X$ means that $B_X(G)\subseteq B'_X(G)$ for every graph $G$. 
The notation $\B_X\leq \B'_X$ means that for every graph $G$ and $R'\in B'_X(G)$, there exists $R\in B_X(G)$ such that $R\subseteq R'$.  
\end{defn}

Let $X$ be a super $X$-set parameter. Then $\B_X\subseteq \B_X^{\max}$ and $\B_X^{\min}\le \B_X$ for every $X$-blocking family $\B_X$.  Since the minimal blocking sets for domination and standard zero forcing are minimal neighborhoods and minimal forts, respectively,  $\B_D^{\min}\subseteq \B_D$ and $\B_{\Z}^{\min}\subseteq \B_{\Z}$.

 For any super $X$-set parameter $X$ 
 and any $X$-blocking family $\B$, we can define a super $Y$-set parameter  (which is usually different from $Y_X$).

 \begin{defn} Let $X$ be a super $X$-set parameter and let $\B_X$ be an $X$-blocking family.  Let $G$ be a graph and $S\subseteq V(G)$. For $u\in S$, a set $R\in B_X(G)$ is a private $B_X$-set for $u$ (relative to $S$) provided $S\cap R = \{u\}$. We say that $S$ is a \emph{$B_X$-irredundant set}  of $G$ provided every element of $S$ has a private $B_X$-set  in $G$. Define $Y_{B_X}(G)$ to be the maximum cardinality of a $B_X$-irredundant set  in $G$.
 \end{defn}

\begin{rem}\label{rem: blocking gives sub param}
Let $X$ be a super $X$-set parameter and let $\B_X$ be an $X$-blocking family. Then $\emptyset$ is trivially a $B_X$-irredundant set of $G$ for every graph $G$ (equivalently, $\emptyset$ is a $Y_{B_X}$-set of $G$). Thus, $Y_{B_X}$ is a cohesive parameter.   Furthermore, $Y_{B_X}$ is a sub $Y$-set parameter. 
\end{rem}

 \subsection{Useful irredundance}\label{ss:irred}
   In order for $B_X$-irredundant sets to be useful, more is needed. For example, the empty $X$-blocking family $\B_X^{\emptyset}$ is not very useful. 

    The  properties in the next three definitions  are useful for producing robust $Y_{B_X}$-set parameters from robust $X$-set parameters (see Theorem \ref{p:xirred-robust}). 
 These properties apply to domination irredundance and zero forcing irredundance, and will apply to the other types of irredundance we introduce.

\begin{defn}\label{d:compconsistXblock}
Let $X$ be a component consistent super $X$-set parameter.  An $X$-blocking family $\B_X$  is \emph{component consistent}  if for every graph $G=G_1\du G_2\du \cdots \du G_k$: 
\ben[$(1)$]
\item\label{26c1}  $R_i\in B_X(G_i)$ implies $R_i\in B_X(G)$.
\item\label{26c2}   $R\in B_X(G)$ implies 
($R\cap V(G_i)\in B_X(G_i)$ or $R\cap V(G_i)= \emptyset$ for $i=1,\dots,k$). 
\een\end{defn}

Note that $\B_X^{\min}$ is component consistent for any component consistent super $X$-set parameter $X$, but  $\B_X^{\max}$ usually is not: For example, Let $G_1,G_2\cong P_3$ with  $V(G_1)=\{1,2,3\}$  and $V(G_2)=\{4,5,6\}$ in path order.  Then $R=\{1,2,4\}$ is a domination blocking (since $1$ cannot be dominated), so $R\in \B_\gamma^{\max}$. But $R\cap\{4,5,6\}=\{4\}$ is not a domination blocking set of of $G_2$.
Since blocking sets for domination and standard zero forcing are adjacency based, $\B_D$ and $\B_{\Z}$ are component consistent.

\begin{defn}\label{d:inclusive}
Let $X$ be a  component consistent super $X$-set parameter. A  component consistent  $X$-blocking family $\B_X$  is \emph{inclusive}  if for every connected graph $G$ of order at least two and $v\in V(G)$, there is some blocking set $R\in B_X(G)$ such that $v\in R$. 
\end{defn}

Since $B_D(G)=\{N[v]:v\in V(G)\}$ and $v\in N[v]$, $\B_D$ is  inclusive. For standard zero forcing, $V(G)$ is a fort of $G$, so $\B_{\Z}$ is inclusive.  

   The next theorem is essential to the study of zero forcing irredundance.

\begin{thm} \label{fortthm}  {\rm \cite{BFH19}}
For a graph $G$, $B\subseteq V(G)$ is a zero forcing set if and only if $B\cap F\ne \emptyset$ for every fort of $G$.
\end{thm} 

Observe that the domination parallel of Theorem \ref{fortthm} is immediate:  For a graph $G$, $U\subseteq V(G)$ is a dominating set if and only if $U$ intersects the closed neighborhood of every vertex.

\begin{defn}\label{d:Xirred} Let $X$ be a super $X$-set parameter and let $ \B_X$ be an $X$-blocking family. Then $ \B_X$ is an \emph{$X$-irredundant blocking family} if for every graph $G$,
\beq\label{d-eq:irred} S\mbox{ is an $X$-set }\Leftrightarrow S\cap R\ne \emptyset\mbox{ for every }R\in B_X(G).\eeq 
Let  $ \B_X$ be an $X$-irredundant blocking family.  Then the $B_X$-irredundant sets are called \emph{$X$Ir-sets}  (or $X$Ir-sets for $\B_X$ if necessary for clarity).  The associated $Y$-set parameter is denoted by $Y_{\B_X}$ or $\XIR$, i.e.,  $\XIR(G)=\max\{|S|: S\mbox{ is a maximal XIr-set}\}$ and called the \emph{upper XIr number}.  The  \emph{lower XIr number} is $\xir(G)=\min\{|S|: S\mbox{ is a maximal XIr-set}\}$. 
\end{defn}

 The next  three results are well known for domination irredundance.  Similar arguments are used here and were used in \cite{ZIR} for zero forcing irredundance.

\begin{rem}\label{r:X+XIr} Let $X$ be a  super $X$-set parameter and let $ \B_X$ be an $X$-irredundant blocking family.  
Let $S\subseteq V(G)$   be both an $X$-set and an XIr-set. Then  $S$ is a minimal $X$-set and a maximal XIr-set: Since every element $x\in S$ has a private blocking set $B_x$, $(S \setminus \{x\})\cap B_x=\emptyset$ and $S \setminus \{x\}$ is not an $X$-set. Since a superset of an $X$-set is a $X$-set, and no proper superset of $S$ is a miminal $X$-set, $S$ is a maximal ZIr-set.
\end{rem}

\begin{prop}\label{p:Xmin-XIrmax}
    Let  $ \B_X$ be an $X$-irredundant blocking family for a  super $X$-set parameter $X$ and let $G$ be a graph. If $S\subseteq V(G)$ is a minimal $X$-set of $G$, then  $S$ is a maximal XIr-set of $G$.
\end{prop}
\bpf 
Suppose that $S$ is a minimal $X$-set. Let $u\in S$. Since $S$ is minimal, $S \setminus \{u\}$ is not an $X$-set. Since $\B_X$ is an $X$-irredundant blocking family, there exists a set $R_u\in \B_X(G)$ such that $(S \setminus \{u\})\cap R_u=\emptyset$.  Since $S$ is an $X$-set,  $S\cap R_u =\{u\}$. Thus $R_u$ is a private $B_X$-set of $u$ relative to $S$. This argument holds for every vertex in $S$ and so $S$ is a XIr-set. 
Furthermore, $S\cap R\ne\emptyset$ for every $\B_X$-set $R$ of $G$ and hence $S$ is a maximal ZIr-set.
\epf

\begin{cor}\label{c:param-interlace}
    Let  $ \B_X$ be an $X$-irredundant blocking family for a component consistent super $X$-set parameter $X$.  Then for every graph $G$,
    \[ \xir(G)\le X(G)\le \xbar(G)\le \XIR(G).\]
\end{cor}

The converse of Proposition \ref{p:Xmin-XIrmax} is not true for many super $X$-set parameters; examples for domination and zero forcing are presented in  \cite{MR21} and \cite{ZIR} (and also in Example \ref{ex:star-Ztar-ZIRtar} for ZIR).

As noted in the the next result, one direction in equation \eqref{d-eq:irred} is true for every $X$-blocking family.

\begin{prop}\label{p:xirred-min}
    Let $X$  be a super $X$-set parameter, let $\B_X$ be an $X$-blocking family, let $G$ be a graph, and let $S\subseteq V(G)$.  

    \ben[$(1)$]
    \item\label{c1p212}   If $S$ is an $X$-set, then $S\cap R\ne \emptyset$ for every $R\in B_X(G)$.
    \item $\B_X$ is an $X$-irredundant blocking family if and only if $\B^{\min}_X\subseteq \B_X$.
    \een
\end{prop}
\bpf 
Let $R$ be an $X$-blocking set.
If $S\cap R=\emptyset$, then $S$ is not an $X$-set because $V(G)\setminus R$ is not an $X$-set and $S\subseteq V(G)\setminus R$. 

 Suppose first that $\B_X$ is $X$-irredundant.  We show that for every graph $G$, $B^{\min}_X(G)\subseteq B_X(G)$.
Let $R_0\in B^{\min}_X(G)$ 
and define $S_0=V(G)\setminus R_0$, so $S_0$ is not an $X$-set. Thus there is some $R\in \B_X$ such that  $S_0\cap R=\emptyset$.  This implies $R\subseteq R_0$.  Since $R_0$ is a minimal $X$-blocking set, $R_0=R\in B_X(G)$. 

Now suppose $\B^{\min}_X\subseteq \B_X$.  By \eqref{c1p212}, to show $\B_X$ is an $X$-irredundant blocking family it suffices to show that for every graph $G$ and set $S\subseteq V(G)$, $S\cap R\ne \emptyset$ for every  $R\in B_X(G)$ implies $S$ is an $X$-set of $G$. To establish the contrapositive, assume $S$ is not an $X$-set of $G$, so $V\setminus S$ is an $X$-blocking set (but not necessarily in $B_X(G)$). Then $V\setminus S$ contains some minimal $X$-blocking set $R_0$,  $S\cap R_0= \emptyset$, and $R_0\in B^{\min}_X(G)\subseteq B_X(G)$.
\epf

The next corollary is immediate.
\begin{cor}
    Let $X$ be a super $X$-set parameter, let $ \B_X$ be an $X$-irredundant blocking family. If $\B'_X$ is an $X$-blocking family such that $\B_X\subseteq \B'_X$, then $\B'$ is an $X$-irredundant blocking family. 
\end{cor}


\subsection{Isomorphisms of TAR graphs of robust irredundance parameters}\label{ss:TAR}

In this section we recall the definition of the TAR reconfiguration graph for a cohesive vertex set parameter 
and show that $Y_{\B_X}=\XIR$ is a robust $Y$-set parameter for a robust $X$-set parameter with a inclusive component consistent $X$-irredundant blocking family $\B_X$.  This allows automatic application of robust $Y$-set parameter results in \cite{XTAR-YTAR} to the TAR graphs of such irredundance parameters.

  For a cohesive $W$-set property $W$,   the {\em token addition and removal  reconfiguration graph (TAR graph)}  of a  base graph  $G$ is the graph defined as follows \cite{XTAR-YTAR}:  The vertex set of the TAR graph is  the set of all  $W$-sets of $G$.  There is an edge between  two vertices  $S_1$ and $S_2$ of the TAR graph of $G$ if and only if  $S_2$ can be obtained from $S_1$ by the addition or removal of exactly one vertex.
 When $X$ is a super $X$-set parameter, the {\em $X$-TAR graph}  of a  base graph  $G$ is denoted by  $\xtar(G)$ and  when $Y$ is a sub $Y$-set parameter, the {\em $Y$-TAR graph}  of a  base graph  $G$ is denoted by  $\ytar(G)$.

 When a specific component consistent family $\B_X$ of  $X$-irredundant blocking sets is understood, the $X$-irredundance TAR graph of $G$ can be denoted by $\xirtar(G)$.
 The next example compares   $\ztar(G)$ and $\zirtar(G)$ for a specific graph $G$.

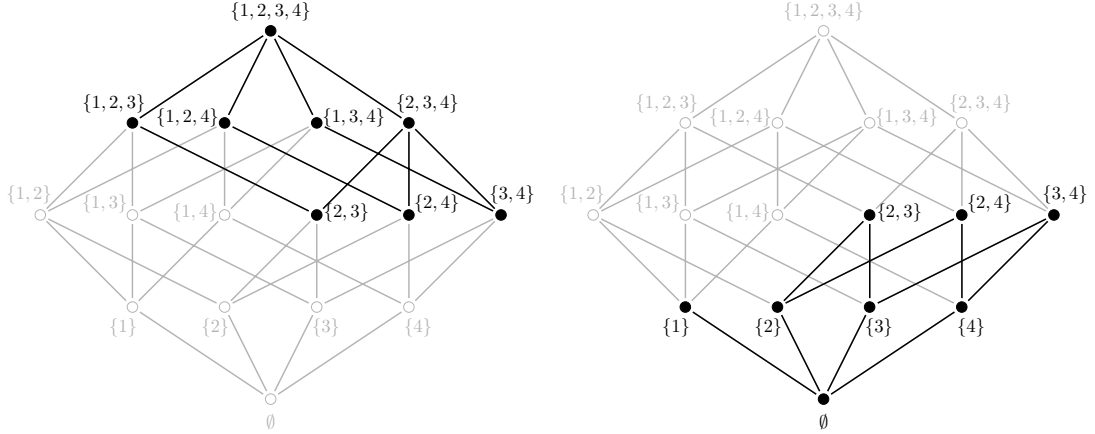
\begin{figure}[!h]
\centering
\centering
\scalebox{0.7}{
\begin{tikzpicture}[scale=1.75,every node/.style={draw,shape=circle,outer sep=2pt,inner sep=1pt,minimum size=.2cm}]		
\node[fill=black, label={[yshift=-20pt]$\{1,2,3,4\}$}]  (1) at (0,2) {};
\node[fill=black, label={[yshift=-15pt, xshift=-10pt]$\{1,2,3\}$}]  (2) at (-1.5,1) {};
\node[fill=black, label={[yshift=-23pt, xshift=-20pt]$\{1,2,4\}$}]  (3) at (-0.5,1) {};
\node[fill=black, label={[yshift=-23pt, xshift=20pt]$\{1,3,4\}$}]  (4) at (0.5,1) {};
\node[fill=black, label={[yshift=-15pt, xshift=10pt]$\{2,3,4\}$}]  (5) at (1.5,1) {};
\node[fill=none, opacity = 0.3, label={[yshift=-10pt, xshift=-6pt, opacity = 0.3]$\{1,2\}$}]  (6) at (-2.5,0) {};
\node[fill=none, opacity = 0.3, label={[yshift=-15pt, xshift=-15pt, opacity = 0.3]$\{1,3\}$}]  (7) at (-1.5,0) {};
\node[fill=black, label={[yshift=-20pt, xshift=16pt]$\{2,3\}$}]  (9) at (0.5,0) {};
\node[fill=none, opacity = 0.3, label={[yshift=-20pt, xshift=-16pt, opacity = 0.3]$\{1,4\}$}]  (8) at (-0.5,0) {};
\node[fill=black, label={[yshift=-15pt, xshift=15pt]$\{2,4\}$}]  (10) at (1.5,0) {};
\node[fill=black, label={[yshift=-10pt, xshift=6pt]$\{3,4\}$}]  (11) at (2.5,0) {};
\node[fill=none, opacity = 0.3, label={[yshift=-28pt, xshift=-5pt, opacity = 0.3]$\{1\}$}]  (12) at (-1.5,-1) {};
\node[fill=none, opacity = 0.3, label={[yshift=-28pt, xshift=-5pt, opacity = 0.3]$\{2\}$}]  (13) at (-0.5,-1) {};
\node[fill=none, opacity = 0.3, label={[yshift=-28pt, xshift=5pt, opacity = 0.3]$\{3\}$}]  (14) at (0.5,-1) {};
\node[fill=none, opacity = 0.3, label={[yshift=-28pt, xshift=5pt, opacity = 0.3]$\{4\}$}]  (15) at (1.5,-1) {};
\node[fill=none, opacity = 0.3, label={[yshift=-25pt, opacity = 0.3]$\emptyset$}]  (16) at (0,-2) {};
		
\draw[thick] 
(1)--(2)--(9)--(5)--(1)--(3)--(10)--(5)--(11)--(4)--(1);
\draw[thick, opacity=0.3]
(1)--(2)--(6)--(12)--(16)--(15)--(11)--(5)--(1)--(3)--(6)--(13)--(16)--(14)--(11)--(4)--(8)--(15)--(10)--(5)--(9)--(2)--(7)--(12)--(8)--(3)--(10)--(13)--(9)--(14)--(7)--(4)--(1);

\node[fill=none, opacity = 0.3, label={[yshift=-20pt, opacity = 0.3]$\{1,2,3,4\}$}]  (1) at (6,2) {};
\node[fill=none, opacity = 0.3, label={[yshift=-15pt, xshift=-10pt, opacity = 0.3]$\{1,2,3\}$}]  (2) at (4.5,1) {};
\node[fill=none, opacity = 0.3, label={[yshift=-23pt, xshift=-20pt, opacity = 0.3]$\{1,2,4\}$}]  (3) at (5.5,1) {};
\node[fill=none, opacity = 0.3, label={[yshift=-23pt, xshift=20pt, opacity = 0.3]$\{1,3,4\}$}]  (4) at (6.5,1) {};
\node[fill=none, opacity = 0.3, label={[yshift=-15pt, xshift=10pt, opacity = 0.3]$\{2,3,4\}$}]  (5) at (7.5,1) {};
\node[fill=none, opacity = 0.3, label={[yshift=-10pt, xshift=-6pt, opacity = 0.3]$\{1,2\}$}]  (6) at (3.5,0) {};
\node[fill=none, opacity = 0.3, label={[yshift=-15pt, xshift=-15pt, opacity = 0.3]$\{1,3\}$}]  (7) at (4.5,0) {};
\node[fill=black, label={[yshift=-20pt, xshift=16pt]$\{2,3\}$}]  (9) at (6.5,0) {};
\node[fill=none, opacity = 0.3, label={[yshift=-20pt, xshift=-16pt, opacity = 0.3]$\{1,4\}$}]  (8) at (5.5,0) {};
\node[fill=black, label={[yshift=-15pt, xshift=15pt]$\{2,4\}$}]  (10) at (7.5,0) {};
\node[fill=black, label={[yshift=-10pt, xshift=6pt]$\{3,4\}$}]  (11) at (8.5,0) {};
\node[fill=black, label={[yshift=-28pt, xshift=-5pt]$\{1\}$}]  (12) at (4.5,-1) {};
\node[fill=black, label={[yshift=-28pt, xshift=-5pt]$\{2\}$}]  (13) at (5.5,-1) {};
\node[fill=black, label={[yshift=-28pt, xshift=5pt]$\{3\}$}]  (14) at (6.5,-1) {};
\node[fill=black, label={[yshift=-28pt, xshift=5pt]$\{4\}$}]  (15) at (7.5,-1) {};
\node[fill=black, label={[yshift=-25pt]$\emptyset$}]  (16) at (6,-2) {};
		
\draw[thick] 
(12)--(16)--(13)--(9)--(14)--(16)--(15)--(11)--(14);
\draw[thick]
(15)--(10)--(13);
\draw[thick, opacity=0.3]
(1)--(2)--(6)--(12)--(16)--(15)--(11)--(5)--(1)--(3)--(6)--(13)--(16)--(14)--(11)--(4)--(8)--(15)--(10)--(5)--(9)--(2)--(7)--(12)--(8)--(3)--(10)--(13)--(9)--(14)--(7)--(4)--(1);
\end{tikzpicture}}
\caption{The TAR reconfiguration graphs  $\ztar(K_{1,3})$ and $\zirtar(K_{1,3})$ (the light gray parts of each hypercube are not included in the TAR graph)}
\label{f:star-Z-ZIR}
\end{figure}

\begin{ex}\label{ex:star-Ztar-ZIRtar}
  Recall that the minimal zero forcing sets of  $K_{1,3}$ (where $V(K_{1,3})=\{1,2,3,4\}$ and  1 is the vertex of degree three) are $\{2,3\}, \{2,4\},$ and $ \{3,4\}$. 
  The forts of $K_{1,3}$ are $\{2,3\}$, $\{2,4\}$, $\{3,4\}$,  $\{2,3,4\}$, and $\{1, 2,3, 4\}$.    The maximal ZIr-sets of $K_{1,3}$ are $\{1\}, \{2,3\}, \{2,4\},$ and $\{3,4\}$. Observe that $\{1\}$ is a maximal $\Z$Ir-set but not a zero forcing set. Figure \ref{f:star-Z-ZIR} shows  $\ztar(K_{1,3})$ and $\zirtar(K_{1,3})$.  
\end{ex}

Next we show that $X$-irredundance inherits important properties from $X$.

 \begin{lem}\label{l:xirred-compconsist}
Let $X$ be a robust   $X$-set parameter  and fix  a component consistent $X$-irredundant blocking family $\B_X$.  Then $\XIR$, the $X$-irredundance number for $\B_X$,  is a sub $Y$-set parameter that satisfies the Component Consistency axiom.
\end{lem}
\bpf  
As noted in Remark \ref{rem: blocking gives sub param}, $\XIR$ is a sub $Y$-set parameter. Let $G$ be a graph with connected components $G_1, \dots, G_k$.  Since $\B_X$ is component consistent, $R\in B_{X}(G)$ nontrivially intersects $V(G_i)$ if and only if $R\cap G_i\in B_X(G_i)$. 
Suppose that $S$ is a XIr-set of $G$. Then each $v\in S$ has a private $B_X(G)$-set $S_v$ and hence $S_v \cap G_i$ is a private $B_X(G_i)$-set, where $G_i$ is the component of $G$ containing $v$. Thus $S\cap V(G_i)$ is a XIr-set of $G_i$.
Now suppose $S'\subseteq V(G)$ such that $S'\cap V(G_i)$ is a XIr-set for each $G_i$. Since  every $v\in S'$ has a private $B_X(G_i)$-set, which is a $B_X(G)$-set, $S'$ is a XIr-set of $G$.
 \epf

\begin{thm}\label{p:xirred-robust}
Let $X$ be a robust  $X$-set parameter with  an inclusive component consistent $X$-irredundant blocking family $\B_X$.  Then $\XIR$, the $X$-irredundance number for $\B_X$,   is a robust $Y$-set parameter.
\end{thm}
\bpf  
 As established in Lemma \ref{l:xirred-compconsist}, $\XIR$ is a sub $Y$-set parameter that satisfies the Component Consistency axiom. Suppose that $G$ is a connected graph  of order at least two. Since $\B_X$ is  inclusive, for every $v\in V(G)$ there is $R\in B_X(G)$ such that $v\in R$. Hence $\{v\}$ is a XIr-set for every $v\in V(G)$. 
 \epf

\begin{thm}\label{t:main-comp} {\rm\cite{XTAR-YTAR}}
Let $Y$ be a robust $Y$-set parameter and let $G$ and $G'$ be base graphs such $\ytar(G)\cong\ytar(G')$. If $Y(K_1)=1$ or $G$ and $G'$ have no isolated vertices, then $G$ and $G'$ have the same order and there is a relabeling of the vertices of $G'$ such that $G$ and $G'$ have exactly the same $Y$-sets.  
\end{thm}

  The next result is an immediate consequence of Theorem \ref{p:xirred-robust} and Theorem \ref{t:main-comp}.

\begin{thm}\label{t:xirred-main}
Let $X$ be a robust  $X$-set parameter with  an inclusive component consistent $X$-irredundant blocking family $\B_X$ and let $\xirtar(G)$ denote the $X$-irredundance TAR graph for $\B_X$. Let $G$ and $G'$ be  base graphs such  that $\xirtar(G)\cong\xirtar(G')$. If $\XIR(K_1)=1$ or $G$ and $G'$ have no isolated vertices, then $G$ and $G'$ have the same order and there is a relabeling of the vertices of $G'$ such that $G$ and $G'$ have exactly the same XIr-sets for $\B_X$.  
\end{thm}

 As noted in \cite{XTAR-YTAR},  $\DIR$ and $\ZIR$ are robust $Y$-set parameters and Theorem \ref{t:xirred-main} applies to $\DIR$ and $\ZIR$.
 
 The hypothesis that $\B_X$ is inclusive is necessary, as noted in the next remark.

  \begin{rem}
    Let $X$ be a super $X$-set parameter and let $G$ be a graph  of order at least two such that $X(G)=0$. Then there are no $X$-blocking sets for $G$ and $B_X(G)=\emptyset$.  When 
$B_X(G)=\emptyset$, the only $\B_X(G)$ irredundant set is $\emptyset$.  Thus $Y_{\B_X}(G)=0$ and $\ytar_{\B_X}(G)\cong K_1$.  Thus $X$-irredundance is not a robust $Y$-set parameter even if $X$ is robust, as the case for the skew forcing (see Section \ref{ss:skew-ired}).  In this case, not only is there no isomorphism theorem for TAR graphs,  the order of the base graph $G$ cannot  be recovered from $\ytar_{\B_X}(G)$. 
 \end{rem}


\section{$X$-irredundant families of blocking sets  for additional specific parameters}\label{ss:block1x}

In this section we describe the natural blocking sets for  PSD forcing,  skew forcing, and vertex covering, and show that these give component consistent $X$-irredundant families of blocking sets. 
We show that the PSD  forcing and vertex covering irredundance parameters are robust, but skew irredundnace is not.

\subsection{PSD zero forcing}\label{ss:PSD-irred}

 Just as (standard) zero forcing naturally uses a (standard) fort as a blocking set, PSD forcing naturally uses a PSD fort as a blocking set.  PSD forts were defined by Smith, Mikesell, and Hicks in \cite{SMH19_PSDforts}: Let $F\subseteq V(G)$ be nonempty and partition $F=F_1\du\dots\du F_k$ where $G[F_i]$ are the connected components of $G[F]$.  Then $F$  is a \emph{PSD fort} of $G$ if and only if $F_i$ is a (standard) fort of $G$ for $i=1,\dots, k$.  Furthermore, they established the result necessary for using these forts for $\Zp$-irredundance:

\begin{thm} \label{fortthmZp}  {\rm \cite{SMH19_PSDforts}}
For a graph $G$, $S\subseteq V(G)$ is a PSD forcing set if and only if $S\cap F\ne \emptyset$ for every PSD fort of $G$.
\end{thm} 

\begin{defn}
    For a graph $G$, let $B_{\Zp}(G)=\{ F: F\mbox{ is a PSD fort of }G\}$ and $\B_{\Zp}=\{B_{\Zp}(G): G\mbox{ is a graph}\}$. A set $S\subseteq V(G)$ is a \emph{$\Zp$-irredundant set} or \emph{$\Zp$Ir-set} provided every element of $S$ has a private PSD fort. 
 The  \emph{upper} $\Zp$Ir \emph{number} is 
\[\ZpIR(G) = \max\{|S|: S \text{ is a maximal $\Zp$Ir-set}\}\] and the \emph{lower} $\Zp$Ir \emph{number} is 
\[\zpir(G) = \min\{|S| :  \text{ is a maximal $\Zp$Ir-set}\}.\] 
\end{defn}

 A PSD fort $F$ of a graph $G$  is \emph{connected} if $G[F]$ is connected.
\begin{obs}\label{r:PSD-fort} 
Since any PSD fort is a disjoint union of connected PSD forts, for the purpose of irredundance we need consider only connected PSD forts.  That is,   A set $S\subseteq V(G)$ is a \emph{$\Zp$-irredundant set} or \emph{$\Zp$Ir-set} provided every element of $S$ has a private connected PSD fort.
\end{obs}

\begin{rem}\label{r:not-n} 
Suppose $G$ is a graph of order $n\ge 2$ that has an edge. Then  $\zpir(G), \ZpIR(G)\le n-1$ because every PSD fort in a connected component containing an edge must contain at least two vertices. 
\end{rem}

 We first present some examples of determining $\ZpIR(G)$ and $\zpir(G)$ for specific graph families $G$.

 \begin{ex}\label{ex:PSDcomplete}
    Let $K_n$ be a complete graph with vertex set $\{v_1,\dots,v_n\}$.  Since $\Zp(K_n)=\zpbar(K_n)=n-1$, $\ZpIR(K_n)=n-1$.  For $k=1,\dots,n$, the set $S_k=\{v_1,\dots,v_{k-1},v_{k+1}\dots,v_n\}$ is a $\Zp$Ir-set (with private PSD fort $\{v_i, v_k\}$ for $v_i$), so $\zpir(K_n)=n-1$.
\end{ex}

 \begin{ex}\label{Zp-empty}
    Let $n\ge 1$ and recall that $\ol{K_n}$ denotes the empty graph.  Then $\zpir(\ol{K_n})=\ZpIR(\ol{K_n})=n$   (it is well known that $\Zp(\ol{K_n})=\zpbar(\ol{K_n})=n$).
\end{ex}

 \begin{ex}\label{ex:PSDtree}
    Let $T$ be a tree.  Since any nonempty set of vertices is a PSD forcing set, $V(T)$ is the  only PSD fort.  Thus $\zpir(T)=\ZpIR(T)=1$ (it is well known that $\Zp(T)=\zpbar(T)=1$).  
\end{ex}

\begin{ex}\label{Zp-cycle}
    Let $n\ge 4$.  It is well known that for the cycle $C_n$, $\Zp({C_n})=2$ and immediate that $\zpbar({C_n})=2$ since any set of two vertices is a PSD forcing set. Note that for any vertex $v$, $F_{\ol v}=V(C_n)\setminus\{v\}$ is PSD fort of $C_n$.   Thus any set of two vertices is a $\Zp$Ir-set.   Let $F$ be a PSD fort of $C_n$ and suppose $v \notin F$ and $N(v)\cap F\ne \emptyset$. Then, $N(v) \in F$. Since $F$ is a PSD fort and $|N(v)|=2$, $N(v)$ is in the same connected component of $F$, implying $F=F_{\ol v}$. If $S$ is a set of 3 or more vertices, no vertex has a private fort, and thus $\zpir(C_n)=\ZpIR({C_n})=2$.
\end{ex}

\begin{ex}\label{ex-star-PSD} 
Consider $K_{p,q}$ with $2\le p\le q$;  let $A$ denote the partite set of $p$ vertices and $B$ denote the partite set of $q$ vertices. Every vertex of a $\Zp$Ir-set must have a private connected fort, and a connected fort requires a vertex in each part.  Furthermore, if $F$ is a connected fort and  $A\not \subseteq F$, then $|B\cap F|\ge 2$, and similarly  if $B\not\subseteq F$. 
Note that $A$ is a $\Zp$-Ir-set  by choosing $F_a=\{a\}\cup B$ as a private fort for $a\in A$, and similarly for $B$. For any $a_1,a_2\in A$ and $b_1,b_2\in B$, $V(K_{p,q})\setminus\{a_1,a_2,b_1,b_2\}$ is a ZIr-set with private connected  fort $\{v,a_1,a_2,b_1,b_2\}$ for $v$. Thus $\zpir(K_{p,q})=p$ and $\ZpIR(K_{p,q})=\max(q,q+p-4)$. It is known that $\Zp(K_{p,q})=p$ \cite{HLS22book} and $\zpbar(K_{p,q})=q$  \cite{XTAR-YTAR}.    %
\end{ex}

Next we examine relationships between the PSD forcing and PSD irredundance parameters.  The next result is immediate from by  Corollary \ref{c:param-interlace}.
\begin{prop}\label{r:Zp-interlace}
   For every graph $G$,  $\zpir(G)\le \Zp(G)\le \zpbar(G)\le \ZpIR(G)$.   
\end{prop}

 Example \ref{ex-star-PSD} shows that $\zpbar(K_{p,q}) < \ZIR(K_{p,q})$ when $p>4$.  Examples \ref{ex:Zp-z+ir}, \ref{ex:ZpbarnotZp} show that the other parameters in the preceding result are distinct. 

\begin{ex}\label{ex:Zp-z+ir}
   Consider the graph $L$ shown in in Figure \ref{f:Zp-z+ir}.   Note that $L$ is obtained from the ladder graph $P_3\Box P_2$ by adjacent twinning of two degree-2 vertices that have a common neighbor, i.e, vertices 4 and 5 are adjacent twins, as are 7 and 8.  For any graph $G$ with adjacent twins  $u$ and $w$, by \cite[Proposition 9.47]{HLS22book}, $\Zp(G)=\Zp(G-w)+1$.  Thus $\Zp(L)=4$ because $\Zp(P_3\Box P_2)=2$.  The set $S=\{2,3,6\}$ is a ZIr-set with private PSD forts  $\{2,1,4,7\}$, $\{3,1,4,7,8\}$, and $\{6,1,4,5,7\}$.  Note that any private PSD fort for 3 with respect to $S$ must contain 1, 7, 8, and similarly any private PSD fort for 6 with respect to $S$ must contain 1, 4, 5. Thus $S$ is a maximal $\Zp$Ir-set and $\zpir(L)\le 3<\Zp(L)$.  (It is not difficult to see that $\zpir(L)= 3$.)

\begin{figure}[!h]
\centering
\scalebox{1}{
\begin{tikzpicture}[scale=1.3,every node/.style={draw,shape=circle,outer sep=2pt,inner sep=1pt,minimum size=.2cm}]		
\node[fill=none]  (4) at (-2.5,-0.5) {$4$};
\node[fill=none]  (3) at (-2,1.5) {$3$};
\node[fill=none]  (5) at (-1.5,0.5) {$5$};
\node[fill=none]  (1) at (0,0) {$1$};
\node[fill=none]  (2) at (0,1.5) {$2$};
\node[fill=none]  (8) at (2.5,0.5) {$8$};
\node[fill=none]  (7) at (1.5,-0.5) {$7$};
\node[fill=none]  (6) at (2,1.5) {$6$};
		
\draw[thick] (2)--(1)--(5)--(4)--(1)--(7)--(8)--(1);
\draw[thick] (4)--(3)--(2)--(6)--(7);
\draw[thick] (3)--(5);
\draw[thick] (6)--(8);
\end{tikzpicture}}
    \caption{\label{f:Zp-z+ir} The graph $L$ that has $\zpir(L)<\Zp(L)$ }  
\end{figure}
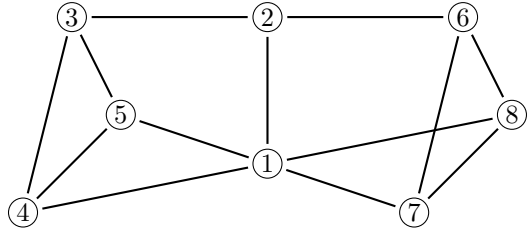

\end{ex}

  \begin{ex}\label{ex:ZpbarnotZp}
   The Brimkov-Carlson graph $B$ shown in Figure \ref{f:BCgraph} was presented   in \cite{BC22} as the smallest member of an infinite family of graphs of graphs $G$ that have   $\Z(G)<\zbar(G)$.  The graph $B$ has the analogous property for PSD forcing: Since $\{1,2,5,11\}$ and  $\{1,2,3,5,10\}$ are both minimal PSD forcing sets,  $\Zp(B)\le 4< 5\le \zpbar(B)$. (It is not difficult to see that $\Zp(B)= 4$ and $ \zpbar(B)=5$.)

\begin{figure}[!h]
\centering
\scalebox{1}{
\begin{tikzpicture}[scale=1.3,every node/.style={draw,shape=circle,outer sep=2pt,inner sep=1pt,minimum size=.51cm}]		
\node[fill=none]  (1) at (0.5,0.87) {$1$};
\node[fill=none]  (2) at (-0.5,0.87) {$2$};
\node[fill=none]  (3) at (-1,0) {$3$};
\node[fill=none]  (4) at (-0.5,-0.87) {$4$};
\node[fill=none]  (5) at (0.5,-0.87) {$5$};
\node[fill=none]  (6) at (1.37,-1.22) {$6$};
\node[fill=none]  (7) at (2.23,-0.87) {$7$};
\node[fill=none]  (8) at (2.59,0) {$8$};
\node[fill=none]  (9) at (2.23,0.87) {$9$};
\node[fill=none]  (10) at (1.37,1.22) {$10$};
\node[fill=none]  (11) at (1.37,0) {$11$};
		
\draw[thick] (11)--(1)--(2)--(3)--(4)--(5)--(11)--(2)--(4)--(11)--(3);
\draw[thick] (11)--(1)--(5)--(11)--(6)--(5);
\draw[thick] (6)--(7)--(11)--(8)--(7);
\draw[thick] (8)--(9)--(11)--(10)--(9);
\draw[thick] (1)--(10);
\draw[thick] (1)--(3)--(5);
\end{tikzpicture}}
   \caption{\label{f:BCgraph} The Brimkov-Carlson graph $G$  has $\Zp(G)<\zpbar(G)$}  
\end{figure}
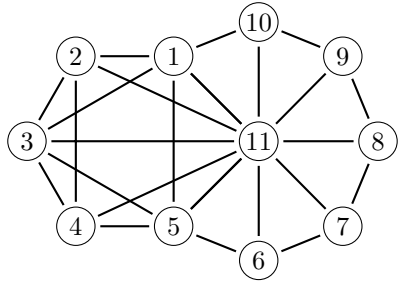

\end{ex}

\begin{rem}\label{r:ZIRandZpIR}
 Observe that any PSD fort is a (standard) fort, but the converse is false. This implies that a $\Zp$Ir-set is a ZIr-set (the converse is false).   
   Thus  $\ZpIR(G)\le \ZIR(G)$.  A path of order $n\ge 5$  shows strict inequality is possible, since $\ZpIR(P_n)=1$ (see Example \ref{ex:PSDtree}) and $\ZIR(P_n)=\lf\frac{n-1}2\rf$ (see  \cite{ZIR}). 
\end{rem}

In contrast to the relationship $\ZpIR(G)\le \ZIR(G)$,  $\zpir(G)$ and $\zir(G)$ are noncomparable, as seen in the next two examples.

\begin{figure}[!h]
\centering
\scalebox{.8}{
\begin{tikzpicture}[scale=2,every node/.style={draw,shape=circle,outer sep=2pt,inner sep=1pt,minimum size=.2cm}]		
\node[fill=none]  (1) at (1.87,-0.5) {$1$};
\node[fill=none]  (2) at (1.87,0.5) {$2$};
\node[fill=none]  (3) at (1,0) {$3$};
\node[fill=none]  (4) at (0,0) {$4$};
\node[fill=none]  (5) at (-0.87,-0.5) {$5$};
\node[fill=none]  (6) at (-0.87,0.5) {$6$};

\draw[thick] (1)--(3)--(4)--(5);
\draw[thick] (4)--(6);
\draw[thick] (2)--(3);
\end{tikzpicture}}
\caption{\label{f:z+ir-zir} A tree $T$ that has $\zpir(T)<\zir(T)$}
\end{figure}
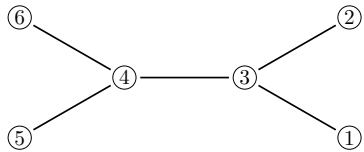
\begin{ex}
    Consider the tree $T$ shown in in Figure \ref{f:z+ir-zir}.  As for every tree, $\zpir(T)=1$. To see that $\zir(T)=2$, note that  
    the complete list of  standard forts of $T$ is 
    $\{1,2\},\{5,6\}, \{1, 2, 3, 5\}, \{1, 2, 3, 6\}, \{1, 2, 5, 6\},$ $\{1, 4, 5, 6\}, \{2, 4, 5, 6\}, \{1, 2, 3, 5, 6\}, \{1, 2, 4, 5, 6\},  $ and $ \{1,2,3,4,5,6\}$ \cite{sage:irred}.  Thus no one vertex is a maximal ZIr-set of $T$, and $\zir(T)\ge 2$.  
 Since $\{3,4\}$ is a maximal ZIr-set of $T$,  $\zir(T)=2$.  
\end{ex}

\begin{ex}    Consider the graph $G$ shown in in Figure \ref{f:z+irmorezir}. We can see that $\zpir(G)=3>2=\zir(G)$:  The complete list of  standard forts of $G$ is 
$\{3, 4\}$,
 $\{5, 6\}$,
 $\{1, 2, 3, 5\}$,
 $\{1, 2, 3, 6\}$,
 $\{1, 2, 4, 5\}$,
 $\{1, 2, 4, 6\}$,
 $\{1, 3, 5, 6\}$,
 $\{1, 4, 5, 6\}$,
 $\{2, 3, 4, 5\}$,
 $\{2, 3, 4, 6\}$,
 $\{3, 4, 5, 6\}$,
 $\{1, 2, 3, 4, 5\}$,
 $\{1, 2, 3, 4, 6\}$,
 $\{1, 2, 3, 5, 6\}$,
 $\{1, 2, 4, 5, 6\}$,
 $\{1, 3, 4, 5, 6\}$,
 $\{2, 3, 4, 5, 6\}$,
 $\{1, 2, 3, 4, 5, 6\}$ \cite{sage:irred}.  Thus $\{1,2\}$ is a maximal ZIr-set for $G$ with private standard forts $\{1, 4, 5, 6\}$ and  $\{2, 3, 4, 5\}$ but these are not PSD forts.  The standard forts that induce connected subgraphs, and are thus connected PSD forts, are 
 $\{3, 4\}$,
 $\{5, 6\}$,
 $\{1, 2, 3, 5\}$,
 $\{1, 2, 3, 6\}$,
 $\{1, 2, 4, 5\}$,
 $\{1, 2, 4, 6\}$,
 $\{1, 2, 3, 4, 5\}$,
 $\{1, 2, 3, 4, 6\}$,
 $\{1, 2, 3, 5, 6\}$,
 $\{1, 2, 4, 5, 6\}$,
 $\{1, 2, 3, 4, 5, 6\}$  ($\{3, 4, 5, 6\}$ is also a PSD fort but is a union of two disjoint PSD forts so is not relevant to finding maximal $\Zp$Ir-sets).  Thus $1$ and $2$ cannot be in the same $\Zp$-Ir-set.  So every   maximal $\Zp$Ir-set for $G$ has at least three elements (one each from $\{1,2\}, \{3,4\}, \{5,6\}$). Thus $\zpir(G)=3$.

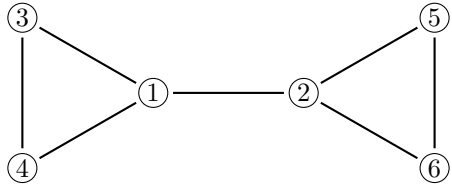
\begin{figure}[!h]
\centering
\scalebox{1}{
\begin{tikzpicture}[scale=2,every node/.style={draw,shape=circle,outer sep=2pt,inner sep=1pt,minimum size=.2cm}]		
\node[fill=none]  (6) at (1.87,-0.5) {$6$};
\node[fill=none]  (5) at (1.87,0.5) {$5$};
\node[fill=none]  (2) at (1,0) {$2$};
\node[fill=none]  (1) at (0,0) {$1$};
\node[fill=none]  (4) at (-0.87,-0.5) {$4$};
\node[fill=none]  (3) at (-0.87,0.5) {$3$};

\draw[thick] (1)--(3)--(4)--(1)--(2)--(5)--(6)--(2);
\end{tikzpicture}}
\caption{\label{f:z+irmorezir}  A graph $G$ that has $\zpir(G)>\zir(G)$} 
\end{figure}

\end{ex}

Next we adapt several results for standard zero forcing irredundance presented in \cite{ZIR}  to PSD forcing.
 Let $G$ be a graph and let  $D\subseteq V(G)$. 
The set $D$  is a \emph{connected $2$-dominating set} if $G[D]$ is connected and for every vertex $v\in V(G)\setminus D$, $v$ is adjacent to at least two distinct vertices in $D$.  The 
\emph{connected $2$-domination number} of $G$ is $\Dctwo(G)=\min\{|D|: D\mbox{ is a connected 2-dominating set}\}$.

\begin{prop}\label{p:2dom-PSD}
Let $G$ be a graph of order $n\ge 3$  with no isolated vertices. 
Then $ n-\Dctwo(G)\le \ZpIR(G)   \leq n-\gamma(G)$.   Both  bounds  sharp.
\end{prop}
\bpf 
The upper bound is immediate from  $\ZpIR(G)\le\ZIR(G)$  and $\ZIR(G)\le n-\gamma(G)$; the latter appears in  \cite{ZIR}.
Let $D$  be a  minimum connected $2$-dominating set and let $S = V(G)\setminus D = \{v_1,\dots,v_k\}$. The set $F_i = D\cup\{v_i\}$ is a private PSD fort of $v_i$ relative to $S$. Therefore $S$ is a $\Zp$Ir-set and $\ZpIR(G)\geq |S|=n-\Dctwo(G)$.
 The graph $K_n$ shows the upper bound is sharp.  For the lower bound, let $G=C_r\vee K_2$,  $V(C_r)=\{v_1,\dots,v_r\}$ with $r\ge 4$, and  $K_2=\{x,y\}$.  For $i=1,\dots,r$, $\{v_i,x,y\}$ is a private PSD fort of  $G$ for $v_i$, so $\ZpIR(G)= r=|V(G)|-2=|V(G)|-\Dctwo(G)$.
\epf

 Next we examine graphs that have extreme values if PSD irredundance numbers. If $G$ does not have an edge, then $G\cong \ol{K_n}$ and $\zpir(\ol{K_n}) = \Zp(\ol{K_n}) = \zpbar(\ol{K_n}) = \ZpIR(\ol{K_n}) = n$. 
It is well known that $\Zp(G),\zpbar(G)\le n-1$ for any graph $G$ of order $n$  that has an edge.

\begin{prop}\label{p:ZpIR_n-1}
Let $G$ be a graph of order $n\ge 2$ with an edge.  
The following are equivalent:
\ben[$(1)$]
\item\label{item:Kn} $G\cong K_{n-r}\du rK_1$ with $n-r\ge 2$. 
\item\label{item:zpir=n-1} $\zpir(G)=n-1$.
\item\label{item:Zp=n-1} $\Zp(G)=n-1$.
\item\label{item:Zpbar=n-1} $\zpbar(G)=n-1$.
\item\label{item:ZpIR=n-1} $\ZpIR(G)= n-1$.
\een
\end{prop}
\bpf 
  The following implications are immediate: 
  \eqref{item:Kn} $\Rightarrow$ \eqref{item:zpir=n-1} $\Rightarrow$ \eqref{item:Zp=n-1} $\Rightarrow$ \eqref{item:Zpbar=n-1}  $\Rightarrow$ \eqref{item:ZpIR=n-1}.  By Remark \ref{r:ZIRandZpIR}, \eqref{item:ZpIR=n-1} implies $\ZIR(G)=n-1$, and  $\ZIR(G)=n-1$ implies   \eqref{item:Kn} by \cite[Proposition 5.2]{ZIR}.
\epf

It is observed in Example \ref{ex:PSDtree}  that if $G$ is a tree, then $\zpir(G)=\Zp(G)=\zpbar(G)=\ZpIR(G)=1$.  Since it is known that $\Zp(G)=1$ if and only if $G$ is a tree, it is immediate that $\ZpIR(G)=1$ (or $\zpbar(G)=1$) implies $G$ is a tree. Next we establish several results which together show that $\zpir(G)=1$ implies $G$ is a tree.
The next lemma and its corollary show that removing a pendent tree from a graph does not change $\zpir$ or $\ZpIR$ (the analogous result is well known for $\Zp$  and $\zpbar$).

\begin{lem}\label{l:treestrip} Let $G$ be a graph, let $x$ be a cut-vertex of $G$,  let   $G[W]$ be a component of $G-x$, and let $\wh W =W\cup\{x\}$, and suppose that $G[\wh W]$ is  a tree. Then:
\ben[$(a)$]
\item\label{l:ts:a} For any PSD fort $F$ of $G$, $F\cap\wh W\ne \emptyset$ implies $\wh W\subseteq F$.
\item\label{l:ts:b} Any $\Zp$Ir-set of $G$ contains at most one element of $\wh W$.
\item\label{l:ts:c} Any maximal $\Zp$Ir-set of $G-W$ is a maximal $\Zp$Ir-set of $G$. For any maximal $\Zp$Ir-set of $G$, there is a maximal $\Zp$Ir-set of $G-W$ of the same size.
\een
\end{lem}

\bpf  
Every vertex of $G[\wh W]$ is a cut-vertex or a leaf of $G$.  
Let $F$ be a PSD fort of $G$ and let $u\in F\cap \wh W$.  If $u$ is a leaf, then its one neighbor is in $F$. So assume $u$ is  a cut-vertex of $G$.  Then every neighbor of $u$ is in a separate component of $G-u$, implying $(N(u)\cap \wh W)\subseteq F$.  Repeating this  argument for each cut-vertex in $F\cap \wh W$ shows $\wh W\subseteq F$, establishing \eqref{l:ts:a}.   
Note that \eqref{l:ts:a} implies \eqref{l:ts:b}, because any private fort of any element in $\wh W$ contains $\wh W$.  

Then \eqref{l:ts:a} implies that any maximal $\Zp$Ir-set of $G-W$ is a maximal $\Zp$Ir-set of $G$.  To complete the proof of \eqref{l:ts:c}, let $S$ be a maximal $\Zp$Ir-set of $G$.  If $W\cap S=\emptyset$, then $S$ is a maximal $\Zp$Ir-set of $G-W$.  If $W\cap S=\{w\}$, then $S\setminus\{w\}\cup\{x\}$ is a maximal $\Zp$Ir-set of $G-W$ using the same private forts except the private fort of $x$ is obtained from the private fort of $w$ by deleting all vertices of $W$. 
\epf

The next result is an immediate consequence of the preceding lemma.

\begin{cor}\label{c:treestrip} Let $G$ be a graph, let $x$ be a cut-vertex of $G$,  let   $G[W]$ be a component of $G-x$, and let $\wh W =W\cup\{x\}$.  If $G[\wh W]$ is  a tree, then $\zpir(G)=\zpir(G-W)$ and  $\ZpIR(G)=\ZpIR(G-W)$.
\end{cor}

It is well known that any graph of order at least two must have  at least two vertices that are not cut-vertices.  This can be seen by choosing a spanning tree of a connected component and noting that a leaf of the spanning tree cannot be a cut-vertex of $G$.  We use this idea in the next proof.

\begin{prop}\label{p:delta2} Let $G$ be a graph $G$ such that $\delta(G)\ge 2$.  Then $\zpir(G)\ge 2$.
\end{prop}
\bpf If  $G$ is not connected, then $\zpir(G)\ge 2$ (since each component contributes to a maximal $\Zp$Ir-set), so assume $G$ is connected. We show that any vertex $u$ of $G$ can be in a $\Zp$Ir-set with  some other vertex, so $\{u\}$ is not maximal $\Zp$Ir-set for any vertex $u$, and thus $\zpir(G)\ge 2$.

Observe that for any vertex  $v$ that is not a cut-vertex, the set $F=V(G)\setminus\{v\}$ is a fort, because $G[F]$ is connected and $v$ has at least two neighbors in $F$ since $\delta(G)\ge 2$. Since $G$ has at least two vertices that are not cut-vertices, for any $u$ that is not a cut-vertex, there is some other $w$ that is not a cut-vertex, and $\{u,w\}$ is a $\Zp$Ir-set.

Now suppose $u$ is a cut-vertex.  Choose a spanning tree of $G$. Let $w$ be a leaf of the spanning tree,  and let $x$ be the closest vertex to $w$ that is a cut-vertex of $G$ and is on the unique path in the spanning tree  from $w$ to $u$.  Let $C=\{c_1=w, c_2, \dots, c_k\}$ be the vertices of the component of $G-x$ that contains $w$ and let $F_w= C$.    Since no vertex on the path from $w$ to $x$ is a cut-vertex of $G$, $x$ has two neighbors in $F_w$.  Thus $F_w$ is a fort, and it  contains $w$ and does not contain $u$. Since $w$ is not a cut-vertex,  the set $F_u=V(G)\setminus\{w\}$ is a fort, and it contains $u$ and does not contain $w$.  Thus $\{u,w\}$ is a $\Zp$Ir-set.
 \epf

\begin{cor}\label{c:zpir-tree} 
Let $G$ be a graph of order $n\ge 2$ with an edge.  
The following are equivalent:
\ben[$(1)$]
\item\label{item:tree} $G$ is a tree. 
\item\label{item:zpir=1} $\zpir(G)=1$.
\item\label{item:Zp=1} $\Zp(G)=1$.
\item\label{item:Zpbar=1} $\zpbar(G)=1$.
\item\label{item:ZpIR=1} $\ZpIR(G)= 1$.
\een
\end{cor}
\bpf It is observed in Example \ref{ex:PSDtree} that $\ZpIR(G)=1$ for every tree $G$.  Furthermore,  $\ZpIR(G)= 1\Rightarrow \zpbar(G)=1\Rightarrow\Zp(G)=1 \Rightarrow \zpir(G)=1$.

To show  that $\zpir(G)=1$ implies $G$ is a tree, 
suppose $G$ is not a tree.  We show this implies $\zpir(G)\ge 2$. If  $G$ is not connected, then $\zpir(G)\ge 2$, so assume $G$ is connected.  Since $G$ is not a tree,  $G$ contains a cycle.  Apply Corollary \ref{c:treestrip} repeatedly to obtain a graph $\wh G$ such that $\zpir(\wh G)= \zpir(G)$ and $\delta(\wh G)\ge 2$.  By Proposition \ref{p:delta2}, this implies $\zpir(G)=\zpir(\wh G)\ge 2$.  
\epf

For all graphs $G$, it is known that $\delta(G)\le \zir(G)$ (see \cite[Corollary 2.11]{ZIR}) and  $\delta(G)\le \Zp(G)$  \cite[Theorem 9.38]{HLS22book}.  We have $\delta(G)\le \zpir(G)$ when $\delta(G)\le 2$ by Proposition \ref{p:delta2} and Corollary \ref{c:zpir-tree}. This leads to the next question.

\begin{quest}\label{q:delta-Z+} Is $\delta(G)\le \zpir(G)$ for all $G$? 
\end{quest}

Finally we consider the TAR graph for PSD irredundance.  Denote the $\Zp$-irredundance TAR graph of $G$ by $\zpirtar(G)$. Since $V(G)$ is a PSD fort for every graph $G$, $\B_{\Zp}$ is inclusive.  Since PSD forts are defined by adjacency,  $\B_{\Zp}$ is component consistent. Thus $\ZpIR$ is a robust $Y$-set parameter. The next result then follows immediately from Theorem \ref{t:xirred-main}.

\begin{thm}\label{t:Zpirred-main}
Let $G$ and $G'$ be connected base graphs that have no isolated vertices. If $\zpirtar(G)\cong\zpirtar(G')$, then $G$ and $G'$ have the same order and there is a relabeling of the vertices of $G'$ such that $G$ and $G'$ have exactly the same $\Zp$Ir-sets.
\end{thm}

\subsection{Skew zero forcing}\label{ss:skew-ired}

Skew forts were defined recently by Bong et al.~in \cite{skewTErecon}, where they also proved the skew version of the theorem needed to use forts for irredundance. Let $G$ be a graph. A nonempty subset  $F\subseteq V(G)$  is a skew fort of $G$ if $|N(v) \cap F|\ne 1$  for each vertex $v\in V(G)$.
 
 \begin{thm} \label{skewfortthm}  {\rm \cite[Theorem 2.2.9]{skewTErecon}}
For a graph $G$, $S\subseteq V(G)$ is a skew forcing set if and only if $S\cap F\ne \emptyset$ for every skew fort of $G$.
\end{thm} 

 \begin{defn}
    For a graph $G$, let $B_{\Zs}(G)=\{ F: F\mbox{ is a skew fort of }G\}$ and $\B_{\Zs}=\{B_{\Zs}(G): G\mbox{ is a graph}\}$. A set $S\subseteq V(G)$ is a \emph{$\Zs$-irredundant set} or \emph{$\Zs$Ir-set} provided every element of $S$ has a private skew fort. 
 The  \emph{upper} $\Zs$Ir \emph{number} is 
\[\ZsIR(G) = \max\{|S|: S \text{ is a maximal $\Zs$Ir-set}\}\] and the \emph{lower} $\Zs$Ir \emph{number} is 
\[\zsir(G) = \min\{|S| :  \text{ is a maximal $\Zs$Ir-set}\}.\] 
\end{defn}

 \begin{obs}
    In a connected graph of order at least two, a skew fort $F$ must be an  set of two or more independent  vertices, or $F$ must have at least three vertices and $G[F]$ must have at least three edges.  If $|F|=2$ then the vertices in $F$ are independent twins. 
\end{obs}

First we determine upper and lower skew irredundance numbers for several well-known families of graphs.

\begin{ex}\label{ex-Kn-skew}
    Consider the complete graph $K_n$.  It is well known that $\Zs(K_n)=n-2$  (see, e.g., \cite{IMA10}) and known that  $\zsbar(K_n)=n-2$ \cite{XTAR-YTAR}.  Since any three vertices form a minimal fort,  $\zsir(K_n)=\ZsIR(K_n)=n-2=\Zs(K_n)=\zsbar(K_n)$.
\end{ex}

\begin{ex}\label{ex:path-skew}
    Consider the path $P_n$.  If $n$ is even then $\zsir(P_n)=\Zs(P_n)=\zsbar(P_n)=\ZsIR(P_n)=0$ by Remark \ref{r:skew0}.  
    Now suppose $n$ is odd,  so $\emptyset$ is not a skew forcing set and $\zsir(P_n)=\Zs(P_n)=1$.  Furthermore, $P_n$ has a unique skew fort $\{v_1,v_3,\dots,v_{2i+1},\dots,v_n\}$ so $\zsbar(P_n)=\ZsIR(P_n)= 1$.
      The values of $\Zs(P_n)$ and $\zsbar(P_n)$ are known (see, e.g., \cite{IMA10} and \cite{XTAR-YTAR}). 
\end{ex}

\begin{ex}\label{ex-Cn-skew}
    Consider the cycle $C_n$.  It is immediate that any skew fort must contain at least every other vertex.  If $n$ is odd, this implies the only fort is the set of all vertices so $\zsir(C_n)=\Zs(C_n)=\zsbar(C_n)=\ZsIR(C_n)=1$. If $n$ is even, this implies there are three forts, the two partite sets and their union, so $\zsir(C_n)=\Zs(C_n)=\zsbar(C_n)=\ZsIR(C_n)=2$.
    The values of $\Zs(C_n)$ and $\zsbar(C_n)$ 
    are known (see, e.g., \cite{IMA10} and \cite{XTAR-YTAR}).  
\end{ex}

 \begin{rem}
    Let $G$ be a bipartite graph with partite sets $A$ and $B$, let $F$ be a skew fort of $G$, and define $F_A=F\cap A$ and $F_B=F\cap B$. If $a\in A$ then $N(a)\cap F_A=\emptyset$ and $|N(a)\cap F_A|=0$.  If $b\in B$, then $N(b)\subseteq A$, so $N(b)\cap F_A=N(b)\cap F$ and $|N(b)\cap F_A|=|N(b)\cap F|\ne 1$.  Thus $F_A$ is a skew fort or empty and similarly for $F_B$.
\end{rem}

\begin{ex}\label{ex-star-skew} 
Consider $K_{p,q}$ with $1\le p\le q$ and $2\le q$; let $A$ denote the partite set of $p$ vertices and $B$ denote the partite set of $q$ vertices.  Since any fort is entirely within one partite set or the union of two such disjoint forts, for $\Zs$Ir-sets we need consider only forts each of which is entirely within one partite set.  In fact, we need only consider sets of two vertices in each partite set (that has at least two vertices).  Thus a set $S$ is a $\Zs$Ir-set if and only if it excludes at least one vertex from each partite set.  Thus $\zsir(K_{p,q})=p+q-2$ (the value $\Zs(K_{p,q})=p+q-2$ is well-known \cite{HLS22book} and implies $\zsbar(K_{p,q})=\ZsIR(K_{p,q})=p+q-2$).
\end{ex}

\begin{rem}
  Since there are examples of connected graphs $G$ of order at least two that have $\ZsIR(G)=0$ (e.g., an even path), $\ZsIR$ is not a robust $Y$-set parameter.  Since $\Zs$ is a robust $X$-set parameter  and $\B_{\Zs}$ is component consistent, both $\ZsIR$ and $\zsir$  satisfy the Component Consistency axiom by Lemma \ref{l:xirred-compconsist}. This implies that if $G$ is the disjoint union of connected components $G_1,\dots,G_k$, then  $\ZsIR(G)=\sum_{i=1}^k\ZsIR(G_i)$ and $\zsir(G)=\sum_{i=1}^k\zsir(G_i)$. 
    (The analogous results are well known for $\Zs$ and immediate for  $\zsbar$ since $\Zs$ is a robust $X$-set parameter \cite{XTAR-YTAR}.) 
\end{rem}

 As defined in \cite{skewTErecon},    the \emph{skew  closure} $\cls(G,S)$ of a set $S\subseteq V(G)$ in a graph $G$ is the set of blue vertices obtained by starting with exactly the vertices in $S$ blue and applying the skew color change rule until no more skew forces are possible.

\begin{prop} {\rm \cite{skewTErecon}} If $G$ is a graph, then $V(G)\setminus \cls(G,\emptyset)$ is the union of all skew forts. 
\end{prop}

 \begin{cor}
    If $G$ is a graph and $S$ is a $\Zs$Ir-set, then $S\subseteq V(G)\setminus \cls(G,\emptyset)$. 
\end{cor}

Next we examine relationships between skew forcing and skew irredundance parameters.  The next result is an immediate consequnce of  Corollary \ref{c:param-interlace}.

\begin{prop}\label{r:Zs-interlace}
   For every graph $G$,  $\zsir(G)\le \Zs(G)\le \zsbar(G)\le \ZsIR(G)$.   
\end{prop}

Examples  \ref{ex:ZsbarnotZ-IR} and \ref{ex:zsir-notZs} and Remark \ref{Zsbar-less-ZsIR} show that all four parameters in the preceding result are distinct. 

\begin{ex}\label{ex:ZsbarnotZ-IR}
     Let $G$ be the graph shown in Figure \ref{f:fan}. Since the minimal skew forcing sets of $G$ are $\{0\}, \{3\},$ and $\{4\}$, $\zsbar(G)=1$.  Since $\{1,2\}$ is a $\Zs$Ir-set with private forts  $\{0, 1, 3, 4\}$  and $\{0, 2, 3, 4\}$, $\ZsIR(G)\ge 2$, so $\zsbar(G)=1<2\le\ZsIR(G)$ (in fact, $\ZsIR(G)= 2$ \cite{sage:irred}).
\end{ex}

\begin{figure}[h!]\centering 

\scalebox{.8}{\begin{tikzpicture}[scale=2,every node/.style={draw,shape=circle,outer sep=2pt,inner sep=1pt,minimum size=.2cm}]		
\node[fill=none]  (0) at (0.5,0.87) {$0$};
\node[fill=none]  (1) at (2,0) {$1$};
\node[fill=none]  (2) at (0,0) {$2$};
\node[fill=none]  (3) at (1.5,.87) {$3$};
\node[fill=none]  (4) at (1,0) {$4$};

\draw[thick] (0)--(2)--(4)--(0)--(3)--(1)--(4)--(3);
\end{tikzpicture}}

\caption{\label{f:fan} A graph $G$ with $\zsbar(G)<\ZsIR(G)$}
\end{figure}

\begin{ex}\label{ex:zsir-notZs}
     Let $G$ be the graph shown in Figure \ref{f:lower-skew}. Observe that vertices $1$ and $2$ are independent twins, so $\Zs(G)=\Zs(G')+1$ where $G'$ is the graph obtained from $G$ by deleting vertex 2 \cite[Proposition 9.87]{HLS22book}.  Since no one vertex of $G'$ is a skew forcing set and $\{0,4\}$ is a skew forcing set of $G'$, $\Zs(G')=2$ and $\Zs(G)=3$.  A complete list of all skew forts of $G$ is
    $ \{1, 2\},
 \{4, 6\},
 \{1, 2, 4, 6\},$
$\{0, 1, 3, 4, 5\},$ 
 $\{0, 1, 3, 5, 6\},
 \{0, 2, 3, 4, 5\},$
$ \{0, 2, 3, 5, 6\},
 \{1, 2, 3, 4, 6\},
 \{0, 1, 2, 3, 4, 5\},$
$ 
 \{0, 1, 2, 3, 5, 6\},$  
 $\{0, 1, 2, 4, 5, 6\},$
$ \{0, 1, 3, 4, 5, 6\},$
$ 
 \{0, 2, 3, 4, 5, 6\},$ $
 \{0,$ $1, 2, 3, 4, 5, 6\}$  \cite{sage:irred}.
Let $S=\{3,5\}$. The only possible private fort for 3 (relative to $S$) is $ \{1, 2, 3, 4, 6\}$ and the only possible private fort for 5 (relative to $S$) is   $\{0, 1, 2, 4, 5, 6\}$, so $S$ is a maximal $\Zs$Ir-set. Thus $\zsir(G)=2<3=\Zs(G)$.
\end{ex}
\begin{figure}[h!]\centering 

\scalebox{.8}{\begin{tikzpicture}[scale=2,every node/.style={draw,shape=circle,outer sep=2pt,inner sep=1pt,minimum size=.2cm}]		
\node[fill=none]  (0) at (0.69,-0.95) {$0$};
\node[fill=none]  (1) at (0.94,0.2) {$1$};
\node[fill=none]  (2) at (1.81,0.59) {$2$};
\node[fill=none]  (3) at (0,0) {$3$};
\node[fill=none]  (4) at (0.94,-0.2) {$4$};
\node[fill=none]  (5) at (0.69,0.95) {$5$};
\node[fill=none]  (6) at (1.81,-0.59) {$6$};

\draw[thick] (1)--(5)--(2)--(6)--(1)--(4)--(2)--(6)--(0)--(3)--(5);
\draw[thick] (0)--(4);
\end{tikzpicture}}

\caption{\label{f:lower-skew} A graph $G$ with $\zsir(G)<\Zs(G)$}
\end{figure}

 Next we characterize graphs having extreme skew irredundance number.

\begin{rem}\label{r:skew0}\label{r:skew-not-robust}
    Let $G$ be a graph such that $\Zs(G)=0$.  Since there are no forts, $\ZsIR(G)=0$ and  $0\le \zsir(G)\le \Zs(G)\le \zsbar(G)\le \ZsIR(G)=0$.    If $\Zs(G')>0$, then $G'$ has at least one fort so $\zsir(G)\ge 1$.  Thus
    \[\emptyset \mbox{ is a skew forcing set }\Leftrightarrow\zsir(G)=0\Leftrightarrow \Zs(G)=0\Leftrightarrow \zsbar(G)=0\Leftrightarrow\ZsIR(G)=0.\]

    It is well known (and easy to see) that $\Zs(\ol{K_n})=n$.  Since an isolated vertex is a skew fort, $\zsir(\ol{K_n})=\ZsIR(\ol{K_n})=n$.
\end{rem}

\begin{rem}\label{r:skew_n-2}
    Let $G$ be a  graph of order $n\ge 2$ that contains an edge.  We can see  that $\ZsIR(G)\le n-2$:   Let $G'$ be a connected component of order $n'\ge 2$ and  note that  it suffices to show that $\ZsIR(G')\le n'-2$ since $\ZsIR$ sums over the connected components.  Suppose to the contrary that 
    $\{v_1,\dots,v_{n'-1}\}$ is $\Zs$Ir-set of $G'$.  Since $|V(G')|=n'$, there is a vertex $w$ such that 
    $\{v_i,w\}$ is a private fort for $v_i$ for every $i=1,\dots,n-1$.  Then $v_i$ is not adjacent to $w$.  But this means $w$ is isolated, a contradiction since $G'$ was assumed to be a connected component.  
\end{rem}

\begin{rem}
For  graph $G$ of order $n\ge 2$ that contains an edge, it is known that $\Zs(G)= n-2$ if and only if $G\cong K_{n_1,\dots,n_t}\du rK_1$ \cite[Theorem 9.76]{HLS22book} (recall $K_n\cong K_{\underbrace{1,\dots,1}_{n \ \mbox{\scriptsize  times}}}$).

 Since a set $S\subseteq V(K_{n_1,\dots,n_t})$ is a $\Zs$Ir-set if and only if it excludes at least one vertex from each of at least two partite sets,  $\zsir(K_{n_1,\dots,n_t})=n-2$ where     
$n=n_1+\dots+n_t$.  Since $\zsir(G)=n-2$ implies $\Zs(G)=n-2$, this implies  $\zsir(G)=n-2$ is equivalent to $G\cong K_{n_1,\dots,n_t}\du rK_1$.
\end{rem}
Theorem \ref{p:skewIr=n-2}  shows neither  $\zsbar(G)=n-2$ nor $\ZsIR(G)=n-2$  implies $G=K_{n_1,\dots,n_t}\du rK_1$.

\begin{thm}\label{p:skewIr=n-2}
    Let $G$ be a connected graph  of order $n\ge 2$.  Then the following are equivalent:
    \ben[$(1)$]
    \item\label{tn-2:c1} $\zsbar(G)=n-2$.
    \item\label{tn-2:c2} $\ZsIR(G)=n-2$.  
    \item\label{tn-2:c3} The vertices of $G$ can be  partitioned as $V(G)=A\du B\du C$ where each of $A$ and $B$ is a nonempty independent set and $G=G[A]\vee G[B]\vee G[C]$.
    \een
\end{thm}
\bpf 

\eqref{tn-2:c3} $\Rightarrow$ \eqref{tn-2:c1}: Let  $G=G[A]\vee G[B]\vee G[C]$ where each of $A$ and $B$ is a nonempty independent set.  
Choose $x\in A$ and $y\in B$, and let $S=V(G)\setminus\{x,y\}$.  
Then $S$ is a skew forcing set because $x$ and $y$ can force each other.  If any vertex in $C$  is removed from $S$, then no forcing can occur.  If a vertex in $A$  is removed from $S$, then $x$ can still force $y$ but no vertex can force $x$, and similarly if a vertex in $B$  is removed from $S$.  Thus $S$ is a minimal skew forcing set.

\eqref{tn-2:c1} $\Rightarrow$ \eqref{tn-2:c2} is immediate from Remarks \ref{r:Zs-interlace} and \ref{r:skew_n-2}.

\eqref{tn-2:c2} $\Rightarrow$ \eqref{tn-2:c3}: Let $S$ be a $\Zs$Ir-set with $|S|=n-2$. Let $\{x,y\}=V(G)\setminus S$. Define $F_x(u)=\{u,x\}$, $F_y(u)=\{u,y\}$, and $F_{xy}(u)=\{u,x,y\}$ for every $u\in S$.  Note that if $F_x(u)$ is a skew fort, then $u$ and $x$ are not adjacent and similarly for $F_y(u)$.  If $F_{xy}(u)$ is a skew fort, then $G[F_{xy}(u)]\cong K_3$.  

Suppose  that $x$ and $y$ are not adjacent.  Then $F_{xy}(u)$ is not a skew fort for any $u\in S$.   Let $z$ be a neighbor of $x$.
Since $x$ and $z$ are not independent, $F_x(z)$ is not a skew fort. But $x$ has only the one neighbor $z$ in $F_y(z)$, so $F_y(z)$ is not a skew fort.  
  So $z$ does not have a private fort, contradicting the assumption that $S$ is a $\Zs$Ir-set.  
  
  So $x$ and $y$ must be adjacent.
  Define $A$ to be the set of $x$ and  all independent twins of $x$, $B$ to be the set of $y$ and  all independent twins of $y$, and $C=V(G)\setminus(A\cup B)$.  We show $G=G[A]\vee G[B]\vee G[C]$. Let $a\in A$ and $b\in B$.  Since $a$ is either $x$ or an independent twin of $x$, $a$ is adjacent to $y$.  Since $b$ is either $y$ or an independent twin of $y$, $b$ is adjacent to $a$.  Thus $G[A\cup B]=G[A]\vee G[B]$.

  Let $c\in C$. If $F_x(c)$ were a fort, then $c$ would be an independent twin of $x$ and thus would be in $A$ rather than $C$.  So $F_x(c)$ is not a skew fort, and similarly $F_y(c)$ is not a skew fort.  Thus the private skew fort of $c$ is $F_{xy}(c)$, which implies $c$ is a neighbor of both $x$ and $y$.  Therefore, $c$ is a neighbor  of every element of $A$ and is a neighbor of  every element of $B$ and $G=G[A]\vee G[B]\vee G[C]$.
 \epf

\begin{cor}
    For a graph $G$ of order $n$, 
    \[\zsbar(G)=n-2\Leftrightarrow\ZsIR(G)=n-2\Leftrightarrow G=(G[A]\vee G[B]\vee G[C])\du rK_1\] where each of $A$ and $B$ is a nonempty independent set.
\end{cor}

\begin{rem}\label{Zsbar-less-ZsIR}
Since $H\vee K_2\cong K_1\vee K_1 \vee H$,    Theorem \ref{p:skewIr=n-2} shows
  $\zsbar(H\vee K_2)=\ZsIR(H\vee {K_2})=|V(H\vee {K_2})|-2=|V(H)|$ for   any  graph $H$. 
  If  $\Zs(H)<|V(H)|-2$ and $H$ has no isolated vertices, then $\Zs(H\vee K_2)<\zsbar(H\vee K_2)$ because  
  $\Zs(G'\vee K_1)= \Zs(G')+1$ for any graph $G'$ with no isolated vertices \cite[Proposition 9.88]{HLS22book}.
\end{rem}

Even though the skew forcing is not robust, there are  results from zero forcing irredundance (and PSD irredundance) that can be adapted to skew irredundance with similar proofs, such as  the next result. 
Let $G$ be a graph.  A set $D\subseteq V(G)$ is a  total 2-dominating set of $G$ if every vertex of $G$ is adjacent to at least two vertices of $D$.  Define  $\DtTwo(G)=\min\{|D|: D \mbox{ is a total 2-dominating set of }G\}$. 

\begin{prop}\label{p:2dom-skew}
Let $G$ be a graph of order $n\ge 3$ with no isolated vertices. 
Then 
$ n-\DtTwo(G)\le \ZsIR(G) \le n-\gamma(G)$  
and both bounds are sharp.
\end{prop}
\bpf 
The upper bound is immediate from  $\ZsIR(G)\le\ZIR(G)$  and $\ZIR(G)\le n-\gamma(G)$; the latter appears in  \cite{ZIR}. The complete bipartite graph $K_{p,q}$ with $2\le p\le q$ shows this bound is sharp.

Let $D$ be a total 2-dominating set of $G$ and let $S = V(G)\setminus D = \{v_1,\dots,v_k\}$. The set $F_i = D\cup\{v_i\}$ is a private skew fort of $v_i$ relative to $S$. 
Therefore $S$ is a $\Zs$Ir-set and $\ZsIR(G)\geq |S|=n-\DtTwo(G)$. 

To see that  the bound is sharp,   construct a graph $G$ from three graphs $H_i, i=1,2,3$ and a $K_3$ with vertices $x_i,i=1,2,3$ as follows: Join $H_i$ to $K_3[\{x_j,x_k\}]$ where $i\ne j,k$.
Figure \ref{f:Dt2} shows the  graph constructed when $H_1=K_2$, $H_2=2K_1$, and $H_3=P_3$.  Any total 2-dominating set necessarily contains at least three vertices and  $\{x_1,x_2,x_3\}$ is a total 2-dominating set for $G$, so $\DtTwo(G)=3$ and $\ZsIR(G)\ge |V(G)|-3$. It follows from Theorem \ref{p:skewIr=n-2} that $\ZsIR(G)< |V(G)|-2$.
\epf

\begin{figure}[h!]\centering 	     

\begin{tikzpicture}[scale=1.5]	
\tikzset{every node/.style={draw,shape=circle,outer sep=2pt,inner sep=1pt,minimum size=.4cm}}
\node[fill=none]  (0) at (0.87,0.5) {$x_2$};
\node[fill=none]  (1) at (0,0) {$x_3$};
\node[fill=none]  (2) at (0.87,-0.5) {$x_1$};
\node[fill=none]  (3) at (-0.5,0.87) {};
\node[fill=none]  (4) at (0.37,1.37) {};
\node[fill=none]  (5) at (0.37,-1.37) {};
\node[fill=none]  (6) at (-0.5,-0.87) {};
\node[fill=none]  (7) at (2.4,0) {};
\node[fill=none]  (8) at (1.82,0.81) {};
\node[fill=none]  (9) at (1.82,-0.81) {};

\tikzset{every node/.style={}}

\node[] at (-0.19,1.33) {$H_1$};
\node[] at (-0.19,-1.33) {$H_2$};
\node[] at (2.8,0) {$H_3$};

\draw[thick] (4)--(3)--(1)--(6)--(2)--(5)--(1)--(2)--(9);
\draw[thick] (3)--(0)--(1)--(4)--(0)--(8)--(2)--(0)--(7)--(8)--(2)--(7)--(9)--(0);
\end{tikzpicture}

\caption{\label{f:Dt2} A graph $G$ with $\ZsIR(G)= n-\DtTwo(G)$}  \end{figure}

\begin{rem}\label{r:ZIRandZsIR}  From the definitions, every skew fort is a standard fort.  This implies that a $\Zs$Ir-set is a ZIr-set (the converse is false).   
   Thus  $\ZsIR(G)\le \ZIR(G)$.  A path $P_n$   with $n\ge 5$  shows strict inequality is possible, since $\ZsIR(P_n)\le 1$ (see Example \ref{ex:path-skew}) and $\ZIR(P_n)=\lf\frac{n-1}2\rf$ (see \cite{ZIR}). 
\end{rem}

In contrast to $\ZsIR(G)\le \ZIR(G)$ for all graphs $G$, $\zsir$ and $\zir$ are noncomparable,  as the path and the star  illustrate: 
As shown in Example  \ref{ex:path-skew} and \cite{ZIR}, $\zsir(P_n)=0<1=\zir(P_n)$ for $n$ even. 
As shown in Example  \ref{ex-star-skew} and \cite{ZIR}, $\zsir(K_{1,n-1})=n-2>1=\zir(K_{1,n-1})$ for $n\ge 3$.

 \begin{rem}  The Leaf Stripping Algorithm \cite[Algorithm 9.79]{HLS22book}   removes leaves and their neighbors from $G$ to produce $\emptyset$ (no vertices) or a graph that has no leaves but has the same skew forcing number as $G$; it produces $\emptyset$ if and only if $\Zs(G)=0$.  When applied to a forest it produces $\emptyset$ or a set of isolated vertices.  Applying the Leaf Stripping Algorithm to a forest  always results in the same cardinality of the set of vertices returned. Thus every minimal skew forcing set has the same size and $\Zs(T)=\zsbar(T)$.
\end{rem}

\begin{quest}\label{q:Zs-forrest}  Are all  parameters 
     $\zsir, \Zs, \zsbar, \ZsIR$ equal for very forest?
\end{quest}

A computer search found $\zsir(T)=\Zs(T)=\zsbar(T)=\ZsIR(T)$ for $|V(T)|\le 11$ \cite{sage:irred}.

\subsection{Vertex covering}\label{ss:VC-irred}

Just as domination naturally uses a closed neighborhood as a blocking set, vertex covering naturally uses an edge as a blocking set.

\begin{defn}
    For a graph $G$, let $B_{VC}(G)=\{ e: e\in E(G)\}$ (recall an edge of $G$ is a 2-element subset of $V(G)$) and $\B_{VC}=\{B_{VC}(G): G\mbox{ is a graph}\}$. A set $S\subseteq V(G)$ is a \emph{VC-irredundant set} or \emph{VCIr-set} provided every element of $S$ has a private edge. 
The \emph{upper VCIR number} is \[\VCIR(G) = \max\{|S|: S \text{ is a maximal VCIr-set}\}\] and the \emph{lower vcir number} is \[\vcir(G) = \min\{|S| : S \text{ is a maximal VCIr-set}\}.\]
\end{defn}

 Recall that the vertex cover  number of $G$ is denoted by $\tau(G)$ and the upper vertex cover  number of $G$ is denoted  by $\vcbar(G)$.  

\begin{rem}\label{VCIr-robust}
      Every edge is a minimal blocking set,  so $\B_{VC}=\B_{VC}^{\min}$ and $\B_{VC}$ is $VC$-irredundant by Proposition \ref{p:xirred-min}.  Thus by Corollary \ref{c:param-interlace},     \[ \vcir(G)\le \vc(G)\le \vcbar(G)\le \VCIR(G).\]
  \end{rem}

 We will see in Corollary \ref{c:vcir=tau} that $\vcir(G)=\vc(G)$.  However, the remaining three parameters are distinct, as seen in Example \ref{e:VCIRbigger}.
 By Proposition \ref{p:Xmin-XIrmax}, every minimal vertex cover is a maximal VCIr-set.  The next example shows that the converse need not hold even when the VCIr-set has size $\vcir(G)$.

\begin{ex}\label{ex:VC-P4}
    Consider the path $P_4$ with  $V(P_4)=\{1,2,3,4\}$ in path order, and the VCIr-set $S=\{1,4\}$. The set $S$ is not a vertex cover since the edge $\{2,3\}$ is not covered, but it is $VC$-irredundant. If we add either middle vertex, not all vertices have a private edge and thus $S$ is a maximal VCIr-set. Note that if vertex $1$ is replaced by vertex $2$, the new set $S'=\{2,4\}$ is a vertex cover and $|S|=|S'|$. Furthermore, $\vcir(P_4)=2=\VCIR(P_4)$.
\end{ex}

The ability to exchange one vertex for its neighbor in the maximal VCIr-set $S$ in Example \ref{ex:VC-P4} (possibly more than once) to create a vertex cover generalizes whenever $|S|=\vcir(G)$. For a graph $G$, a maximal VCIr-set $S$ such that  $|S|=\vcir(G)$ is called a vcir-set.

\begin{prop}\label{vcirtocover} Let $G$ be a graph.  
There is a vcir-set of $G$ that is a  vertex cover of $G$.   \end{prop}
\bpf 
 Suppose $S$ is a vcir-set for $G$ that is not a vertex cover.  We show that there is a  vcir-set $S'$ that has fewer uncovered edges. 
Suppose the edge $xy$ is not covered by $S$.  If for every neighbor $w$ of $x$ with $w \in S$, $w$ has a private edge that is not $wx$, then $x$ could be added to $S$ to obtain a larger VCIr-set, contradicting the maximality of $S$.  

Therefore $x$ has a neighbor $z \in S$ such that $zx$ is the only private edge for $z$ (and this statement is true for every endpoint of every edge of $G$ that is not covered by $S$).  Since $zx$ is the only private edge for $z$, $(N(z)\setminus\{x\})\subseteq S$. 

Suppose there exists other vertices $w_1, w_2, \dots, w_j \in N(x)\setminus\{y, z\}$ such that $w_i \in S$ for $1 \leq i \leq j$ and $w_ix$ is the only private edge for $w_i$. Then $(S\cup\{x\})\setminus\{z, w_1, w_2, \dots, w_j\}$  
is a maximal VCIr-set of smaller cardinality, contradicting the minimality of $S$. Thus no such vertices $w_1, w_2, \dots, w_j \in N(x)\setminus\{y, z\}$ exist. 

Thus, $z$ is the only neighbor $u$ of $x$ having the edge $ux$ as its unique private edge.  Define $S'=(S\cup\{x\})\setminus\{z\}$.  Note that edge $xy$ is now covered and no edge has been uncovered.  Furthermore $S'$ is maximal because every  edge $ab$ that is still not covered has retained the property that $a$ has one neighbor $u$ such that $ua$ is the only private edge for $u$ and similarly for $b$, so adding either $a$ or $b$ to $S'$ results in a set that is not a VCIr-set.

Repeated application of this process produces a desired set $S_o$ that is  a vcir-set. 
\epf

The next result is immediate from  Proposition \ref{vcirtocover}  and  Remark \ref{VCIr-robust}.
\begin{cor}\label{c:vcir=tau}
    For any graph $G$, $\vcir(G)=\vc(G)$.
\end{cor}

 Henceforth we will use the well-known notation $\vc(G)$  for this parameter.
 Next we determine upper and lower vertex cover irredundance numbers for several well-known families of graphs.

\begin{rem}\label{rem:vcirKn}
   It is well known  (and easy to see) that for every graph $G$  of order $n$ that has an edge, $\VCIR(G)\le n-1$ and $\tau(G)=n-1$ if and only if $G\cong K_n$. Thus, $\vcbar(K_n)=\VCIR(K_n)=n-1$.  
\end{rem}

\begin{rem}\label{r:VCIrNoIsolates} 
   The empty graph $\ol{K_n}$ has no edges, so $B_{VC}(\ol{K_n})=\emptyset$ and $\tau(\ol{K_n})=\vcbar(\ol{K_n})=\VCIR(\ol{K_n})=0$. Since isolated vertices cannot have a private edge,  they are not in VCIr-sets. It is natural focus the  study of vertex cover irredundance on the class of graphs with no isolated vertices. 
\end{rem}

The next example shows that  $\VCIR(G)$ may be greater than $\vcbar(G)$, which may in turn be greater than $\tau(G)$.

\begin{ex}\label{e:VCIRbigger}
    The complete bipartite graph $K_{p,q}$ with $p \leq q$ and partite sets $A=\{a_1,a_2,\dots,a_p\}$ and $B=\{b_1,b_2,\dots,b_q\}$ has $\vc(K_{p,q})=p$ and $\vcbar(K_{p,q})=q$,  since $A$ and $B$ are the only two minimal vertex covers.  If 
     $S\subseteq A$ or $S\subseteq B$, then  $S$ is a VCIr-set, and $S=A$ and $S=B$ are maximal VCIr-sets.  If $p\ge 2$, then a set $V(K_{p,q})\setminus\{a_x,b_y\}$ for some $1 \leq x \leq p$ and $1 \leq y \leq q$ is also  a VCIr-set: The private edge for each $a_i \neq a_x$ is $a_ib_y$ and the private edge for each $b_j \neq b_y$ is $b_ja_x$. Such a set $S$ is a maximal VCIr-set since adding one of the missing vertices results in a set with an element without a private edge.   Thus  $\VCIR(K_{p,q})=p+q-2$ for $p\ge 2$ and $\VCIR(K_{1,q})=q$.
\end{ex}

 \begin{rem}\label{p:vcir-one}    Let $G$ be a graph of order $n$ with no isolated vertices. Then $\vc(G)=1$ if and only if $G=K_{1,n-1}$ since a single vertex must be incident to every edge. \end{rem}

\begin{ex}\label{ex:VC-path-cyc}
 Label the vertices of $P_n$ or $C_n$ in  path or cycle order.  A subset $S$ of vertices  is a vertex cover   if and only if it does not omit two  consecutive vertices.  By choosing every other vertex we obtain the   well-known results that $\vc(P_n)=\lf\frac{n}{2}\rf$ and $\vc(C_n)=\lc\frac{n}{2}\rc$.  If three consecutive vertices are included (or of an endpoint of a path and its neighbor are included), then the set $S$ cannot be a VCIr-set.  By choosing two of every three vertices starting with one, skipping one, choosing two, skipping one, etc., we 
 see that $\vcbar(P_n)=\VCIR(P_n)=\lc\frac{2n-1}3\rc$.  This also works for the cycle when $n\equiv 0, 2\mod 3$ but not when $n\equiv 1\mod 3$, where one fewer vertex is allowed, so $\vcbar(C_n)=\VCIR(C_n)=\lc\frac{2n-2}3\rc$
\end{ex}

\begin{prop}\label{p:VCIRmax} 
    Let $G$ be a graph on $n$ vertices. Then  $n-\gamma(G)\le \VCIR(G)$ and this bound is sharp.  Furthermore, $\VCIR(G)=n-1$ if and only if $\gamma(G)=1$. 
\end{prop}

\bpf Let $D$ be a minimum dominating set of $G$ and let $S=V(G)\setminus D$.  Then $S$ is a VCIr-set of $G$ because every vertex in $S$ has an edge with a vertex in $D$. Thus $n-\gamma(G)\le \VCIR(G)$.  The complete bipartite graph shows the bound is sharp.

Conversely, let $S$ be a VCIr-set with $|S|=\VCIR(G)=n-1$. Then every $x \in S$ has a private edge. Let $w \notin S$. Then, for every $x \in S$, $wx$ must be the private edge of $x$ with respect to $S$. Thus $N[w]=V(G)$. 
\epf

  By Proposition \ref{p:VCIRmax},  $\VCIR(G)=n-1$ does not imply $G\cong K_n$.
The next result is immediate from \cite{ZIR} and Propositions \ref{p:2dom-PSD}, \ref{p:2dom-skew}, and \ref{p:VCIRmax}.

\begin{cor}\label{c:CVIRbig} Let $G$ be a  graph  of order $n\ge 2$   with no isolated vertices. Then $\VCIR(G)\ge \ZIR(G), \VCIR(G)\ge \ZpIR(G),$ and $\VCIR(G)\ge \ZsIR(G)$.
\end{cor}

\begin{cor} Let $G$ be a  graph of order $n$ with no isolated vertices. Then $\VCIR(G)\ge \frac n 2$ and this bound is sharp.
\end{cor}
\bpf The bound follows since $\gamma(G)\le \frac n 2$ \cite[Theorem 2.1]{HHS98}.  The comb $P_r\circ K_1$ for $r\ge 2$ shows the bound is sharp, since only one of a leaf and its neighbor can be in a VCIr-set.
\epf

  The relationships between the upper and lower parameters for domination, independence and domination irredundance  are expressed in the  \emph{Domination Chain},  which was  introduced in \cite{CHM78} and discussed in \cite{MR21}: 
\beq\label{eqDomChain} \dir{G} \leq \gamma(G) \leq \underline{\alpha}(G) \leq \alpha(G) \leq \ol{\gamma}(G) \leq \DIR(G).\eeq
  Various aspects of the Domination Chain have been studied. In \cite{CM93} it was shown that the values of these six parameters can be arbitrary (within known constraints) for a connected graph. That is, for any positive integers $k_1\leq k_2\leq k_3\leq k_4\leq k_5\leq k_6$ such that $k_1=1$ implies $k_3=1$, $k_4=1$ implies $k_6=1$, and that $k_2\leq 2k_1-1$, there exists a connected graph $G$ such that $\dir(G)=k_1, \gamma(G)=k_2, \underline{\alpha}(G)=k_3, \alpha(G)=k_4, \dbar(G)=k_5$ and $\DIR(G)=k_6$. The conditions under which these parameters are equal have also been studied (see \cite{MR21} for a survey). The complexity of computing parameters in the Domination Chain is studied in \cite{BBCF20domchain} and the references therein.
 Next we consider to what extent $\VCIR$ and $\vcir=\tau$ fit in  this chain.

\begin{prop}\label{p:VCIrandDir} 
 Let $G$ be a graph with no isolated vertices. If $S \subseteq V(G)$ is a DIr-set, then $S$ is a VCIr-set. Thus $\DIR(G)\le\VCIR(G)$.
\end{prop}

\bpf 
Let $S \subseteq V(G)$ be a DIr-set. Then for every $x \in S$, there exists $w \in V(G)$ such that $N[w] \cap S = \{x\}$. If $w \neq x$, the edge $wx$ is a private edge for $x$ with respect to $S$. If $w=x$, there exists $u \in V(G)$ such that  $u$ is a private neighbor of $w=x$, which implies $ux$ is a private edge for $t$ with respect to $S$.  The last statement then follows.
\epf

This result, together with Corollary \ref{c:CVIRbig},  shows that among the irredundance parameters discussed, $\VCIR$ is the largest.
\begin{rem} \label{r:ExtDomChain}
      Proposition \ref{p:VCIrandDir} adds the parameter $\VCIR(G)$ to the Domination Chain  for a graph $G$ with no isolated vertices  to create the \emph{Extended Domination Chain} for graphs with no isolated vertices:
 \[\dir(G) \leq \gamma(G) \leq \underline{\alpha}(G) \leq \alpha(G) \leq \ol{\gamma}(G) \leq \DIR(G) \leq \VCIR(G).\]
\end{rem}

Note that a maximal VCIr-set need not be a minimal dominating set nor a DIr-set, as seen in the next example.

\begin{ex}\label{ex:VC-paw}
    Let $G$ be the Paw graph shown in Figure \ref{f:paw}, so $V(G)=\{1,2,3,4\}$ with vertex $1$ the leaf and vertices $2,3$ of degree 2.  The set $\{1,2,3\}$ is a maximal VCIr-set and a dominating set, but not a minimal dominating set nor a DIr-set.  Furthermore, $\{4\}$ is a minimal dominating set and  a minimal DIr-set.   It is not, however, a maximal VCIr-set, because $\{4,2\}$ is a VCIr-set.  In fact, $\dir(G)=1<2=\vc(G)$ and $\DIR(G)=2<3=\VCIR(G)$.

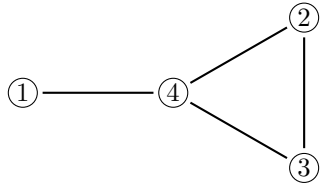
\begin{figure}[!h]
\centering
\scalebox{1}{
\begin{tikzpicture}[scale=2,every node/.style={draw,shape=circle,outer sep=2pt,inner sep=1pt,minimum size=.2cm}]		
\node[fill=none]  (1) at (0,0) {$1$};
\node[fill=none]  (4) at (1,0) {$4$};
\node[fill=none]  (2) at (1.87,0.5) {$2$};
\node[fill=none]  (3) at (1.87,-0.5) {$3$};

\draw[thick] (1)--(4)--(2)--(3)--(4);
\end{tikzpicture}}
\caption{\label{f:paw} The Paw graph} 
\end{figure}

\end{ex}

 \begin{rem} Let $G$ be a graph with no isolated vertices.  Since any vertex cover is a dominating set,  $\tau(G)\ge \gamma(G)$. The Domination Chain then  implies $\dir(G) \leq \gamma(G)\le \tau(G)$.
 \end{rem}

It is  natural  to ask  how $\vc(G)$ compares to the other parameters in the Domination Chain. The next example shows that $\vc(G)$ does not interlace in the Domination Chain because $\vc(G)$ and  the parameters $\ibar(G), \alpha(G), \dbar(G)$, and 
 $\DIR(G)$ are not comparable.

 \begin{ex}
     For the complete graph $K_n$, $\vc(K_n)=n-1$ and $\DIR(K_n)=1$. On the other hand,  the tree $T$ in Figure \ref{f:z+ir-zir} has $\vc(T)=2$  ($\{3,4\}$ is a vertex cover) and  $\ibar(T)=3$ ($\{1,2,4\}$ is a maximal independent set of minimum cardinality). 
 \end{ex}

 Vertex cover is a robust $X$-set parameter. Since every vertex is incident to an edge for every connected graph $G$ of order at least two, $\B_{VC}$ is inclusive.  Since the blocking sets are edges,  $\B_{VC}$ is component consistent. Thus $\VCIR$ is a robust $Y$-set parameter. Given a graph $G$ with no isolated vertices, denote the VCIR-TAR graph of $G$ by $\vcirtar(G)$. Theorem \ref{t:xirred-main} holds for vertex cover irredundance.

\begin{thm}\label{t:vcirmain}
    Let $G$ and $G'$ be base graphs with no isolated vertices such $\vcirtar(G)\cong\vcirtar(G')$. Then $G$ and $G'$ have the same order and there is a relabeling of the vertices of $G'$ such that $G$ and $G'$ have exactly the same VCIr-sets.
\end{thm}

\begin{rem}
    Although $\vcir(G)=\vc(G)$, Example \ref{ex:VC-P4} shows $\vcirtar(P_4)\not\cong\vctar(P_4)$.
\end{rem}

 It was shown in \cite{XTAR-YTAR} that the every graph has a unique VC-TAR graph. This is not the case for vertex cover irredundance.

\begin{ex}
    Let $G=K_2 \sqcup K_2$ with the edges $12$ and $34$. Then $G$ and $P_4$ have the same VCIr-sets. With the vertices of $P_4$ labeled in path order, the maximal VCIr-sets are $\{1,3\}$, $\{1,4\}$, $\{2,3\}$ and $\{2,4\}$ for both graphs. 
\end{ex}

 \begin{rem}
    For graphs with no isolated vertices, every independent set is a Dir-set and every DIr-set is a VCir-set.  Thus the independence TAR graph is an induced subgraph of the domination irredundance TAR graph which is an induced subgraph of the vertex cover irredundance TAR graph. The next example discusses a graph for which all three TAR graphs are different.
    \end{rem}

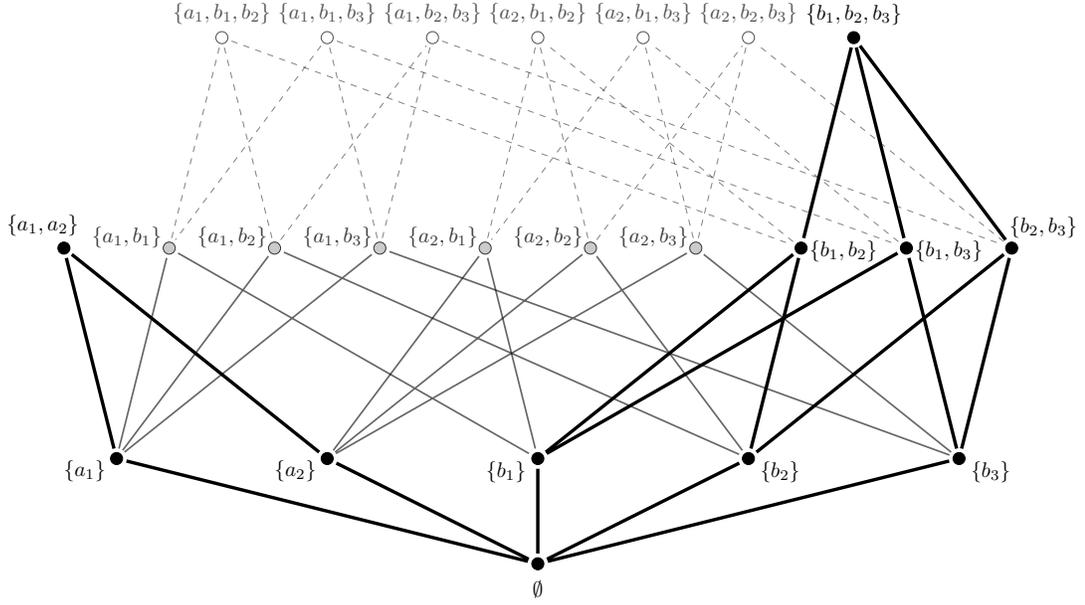
\begin{figure}[!h]
\centering

\scalebox{0.8}{
\begin{tikzpicture}[scale=1.75]

\tikzset{every node/.style={draw,shape=circle,outer sep=2pt,inner sep=1pt,minimum size=.2cm,fill=black}}

\node[label={[yshift=-25pt]$\emptyset$}]  (1) at (0,0) {};
\node[label={[yshift=-25pt, xshift=-15pt]$\{a_1\}$}]  (2) at (-4,1) {};
\node[label={[yshift=-25pt, xshift=-15pt]$\{a_2\}$}]  (3) at (-2,1) {};
\node[label={[yshift=-25pt, xshift=-15pt]$\{b_1\}$}]  (4) at (0,1) {};
\node[label={[yshift=-25pt, xshift=15pt]$\{b_2\}$}]  (5) at (2,1) {};
\node[label={[yshift=-25pt, xshift=15pt]$\{b_3\}$}]  (6) at (4,1) {};
\node[label={[yshift=-15pt, xshift=-10pt]$\{a_1,a_2\}$}]  (7) at (-4.5,3) {};
\node[label={[yshift=-26pt, xshift=20pt]$\{b_1,b_2\}$}]  (8) at (2.5,3) {};
\node[label={[yshift=-26pt, xshift=20pt]$\{b_1,b_3\}$}]  (9) at (3.5,3) {};
\node[label={[yshift=-15pt, xshift=15pt]$\{b_2,b_3\}$}]  (10) at (4.5,3) {};
\node[label={[yshift=-20pt]$\{b_1,b_2,b_3\}$}]  (23) at (3,5) {};

\tikzset{every node/.style={draw,shape=circle,outer sep=2pt,inner sep=1pt,minimum size=.2cm,fill=black!25, opacity=0.8}}

\node[label={[yshift=-20pt, xshift=-20pt]$\{a_1,b_1\}$}]  (11) at (-3.5,3) {};
\node[label={[yshift=-20pt, xshift=-20pt]$\{a_1,b_2\}$}]  (12) at (-2.5,3) {};
\node[label={[yshift=-20pt, xshift=-20pt]$\{a_1,b_3\}$}]  (13) at (-1.5,3) {};
\node[label={[yshift=-20pt, xshift=-20pt]$\{a_2,b_1\}$}]  (14) at (-0.5,3) {};
\node[label={[yshift=-20pt, xshift=-20pt]$\{a_2,b_2\}$}]  (15) at (0.5,3) {};
\node[label={[yshift=-20pt, xshift=-20pt]$\{a_2,b_3\}$}]  (16) at (1.5,3) {};

\tikzset{every node/.style={draw,shape=circle,outer sep=2pt,inner sep=1pt,minimum size=.2cm,fill=none, opacity=0.7}}

\node[label={[yshift=-20pt]$\{a_1,b_1,b_2\}$}]  (17) at (-3,5) {};
\node[label={[yshift=-20pt]$\{a_1,b_1,b_3\}$}]  (18) at (-2,5) {};
\node[label={[yshift=-20pt]$\{a_1,b_2,b_3\}$}]  (19) at (-1,5) {};
\node[label={[yshift=-20pt]$\{a_2,b_1,b_2\}$}]  (20) at (0,5) {};
\node[label={[yshift=-20pt]$\{a_2,b_1,b_3\}$}]  (21) at (1,5) {};
\node[label={[yshift=-20pt]$\{a_2,b_2,b_3\}$}]  (22) at (2,5) {};
		
\draw[ultra thick] 
(1)--(2)--(7)--(3)--(1)--(4)--(8)--(5)--(1)--(6)--(9)--(4);
\draw[ultra thick]
(6)--(10)--(5);
\draw[ultra thick]
(8)--(23)--(9);
\draw[ultra thick]
(10)--(23);

\draw[thick, opacity=0.6]
(4)--(11)--(2)--(13)--(6)--(16)--(3)--(15)--(5)--(12)--(2);
\draw[thick, opacity=0.6]
(3)--(14)--(4);

\draw[dashed, opacity=0.5]
(11)--(17)--(12)--(19)--(13)--(18)--(11);
\draw[dashed, opacity=0.5]
(14)--(20)--(15)--(22)--(16)--(21)--(14);
\draw[dashed, opacity=0.5]
(17)--(8)--(20);
\draw[dashed, opacity=0.5]
(18)--(9)--(21);
\draw[dashed, opacity=0.5]
(19)--(10)--(22);

\end{tikzpicture}
}

\caption{\label{f:K23-VCIR-DIR-TAR}  The $\VCIR, \DIR,$ and $\alpha$ TAR reconfiguration graphs of $K_{2,3}$.  The independence TAR graph is shown in black. The domination irredundance  TAR graph includes both black and gray edges and vertices.  The vertex cover irredundance TAR graph includes all vertices and edges shown.} 
\end{figure}

    \begin{ex}
     Consider the complete bipartite graph $K_{p,q}$ with $p \leq q$ and partite sets $A=\{a_1,a_2,\dots,a_p\}$ and $B=\{b_1,b_2,\dots,b_q\}$. As discussed in Example \ref{e:VCIRbigger}, a set $S\subseteq V(K_{p,q})$ is a VCIr-set if and only if 
     $S\subseteq A$, $S\subseteq B$, or ($A\not\subseteq S$ and $B\not\subseteq S$). It is immediate that $S$ is an independent set if and only if $S\subseteq A$ or $S\subseteq B$.  Furthermore, $S$ is a DIr-set if and only if $S\subseteq A$, $S\subseteq B$, or $S=\{a_i,b_j\}$ for some $,j$ with $1\le i
     \le p, 1\le j\le q$.
     The independence, domination irredundance, and vertex cover irredundance TAR reconfiguration graphs of $K_{2,3}$ are shown in  Figure \ref{f:K23-VCIR-DIR-TAR}.  Observe that all three TAR graphs are different.
    \end{ex}

\section{$X$-irredundant blocking set construction by closure operators}\label{ss:block2}

In this section we introduce a technique for constructing $X$-blocking families that utilizes closure operators.  This technique  can produce component consistent $X$-irredundant blocking families and is  naturally compatible with the propagation processes associated with many super $X$-set parameters  that involve forcing. Let $\mathcal{P}(V)$ denote the power set of $V$.  For any map $\psi:V \to W$ and subset $A\subseteq V$ we write $\psi(A)$ to mean the image of $A$, i.e.,\ $\psi(A) = \{\psi(a) : a\in A\}$.  A  map $\varphi:V \to V$ is a \emph{closure operator} if $A\subseteq \varphi(A)$ for all $A\subseteq V$, $A_1\subseteq A_2$ implies $\varphi(A_1) \subseteq \varphi(A_2)$, and $\varphi \circ \varphi = \varphi$. 

\begin{defn}  
For each graph $G$ assign a  closure operator $\varphi_{G} : \mathcal{P}(V(G)) \to \mathcal{P}(V(G))$. The set $\Phi = \{\varphi_{G} : G \text{ is a graph}\}$ is a \emph{closure family} provided for every isomorphism $\psi: V(G)\to V(H)$, 
$\psi\circ \varphi_{G} = \varphi_{H}$.
\end{defn}

\begin{defn}
    A closure family $\Phi$ is \emph{component consistent} if for every graph $G$, every  component $G_i$ of $G$, and $A\subseteq V(G_i)$, $\varphi_G(A)\subseteq V(G_i)$.
\end{defn}

 \begin{defn}\label{d:inclusive2}
    A component consistent closure family $\Phi$ is \emph{inclusive} if    for every connected graph $G$ of order at least two and $v\in V(G)$, there is some $A\subseteq V(G)$ such that  $v\not\in \varphi(A)$.
\end{defn}

For a component consistent closure family $\Phi$, note that $\varphi_G(\emptyset)=\emptyset$ for every connected graph $G$ of order at least two implies that $\Phi$ is inclusive.

\begin{defn}
Let $X$ be a super $X$-set parameter. For each graph $G$ assign a  closure operator $\varphi_{X,G} : \mathcal{P}(V(G)) \to \mathcal{P}(V(G))$. We refer to $\Phi = \{\varphi_{X,G} : G \text{ is a graph}\}$ as a \emph{closure family} and say that $\Phi$ is $X$-\emph{compliant} provided that  for every graph $G$, 
$S$ is an $X$-set of $G$ if and only if $\varphi_{X,G}(S) = V(G)$.

For an $X$-compliant closure family $\Phi$, define $B_{X,\Phi}(G) = \{ V(G) \setminus \varphi_{X,G}(A) :  A \text{ is not an $X$-set}\}$ and $\B_{X,\Phi} = \{B_{X,\Phi}(G) : G \text{ is a graph}\}$.
\end{defn}

\begin{obs}\label{o:Xcomply}
 The definition of $B_{X,\Phi}(G)$  is equivalent to 
\[B_{X,\Phi}(G) = \{ V(G) \setminus  A:   \mbox{ such that } A=\varphi_{X,G}(A) \mbox{ and $A$ is not an $X$-set}\}.\] 
For an $X$-compliant closure family $\Phi$, $R\in B_{x,\Phi}(G)$ if and only if both $\varphi_{X,G}(V\setminus R)=V\setminus R$ and $R\ne\emptyset$.
\end{obs}

\begin{obs}
If $X$ is a super $X$-set parameter and $\Phi$ is an $X$-compliant closure family, then $\B_{X,\Phi}$ is an $X$-blocking family.    
\end{obs}

\begin{prop}\label{o:clos-compconst}
Let $X$ be a  component consistent super $X$-set parameter and let $\Phi$ be a   $X$-compliant component consistent closure family. Then $\B_{X,\Phi}$ is  component consistent (as an $X$-blocking family). 
\end{prop}
 \bpf Let the components of $G$ be $G_1,\dots, G_k$, let  $V=V(G)$,  and $V_i=V(G_i),i=1,\dots,k$. 
Since  $\Phi$ is component consistent, $\varphi_{X,G_i}={\varphi_{X,G}|}_{V_i}$, where $\psi|_W$ denotes the restriction of the function $\psi$ to the set $W$.  We use $\varphi$ to denote $\varphi_{X,G}$ or $\varphi_{X,G_i}$ as needed.

To establish condition \eqref{26c1} of Definition \ref{d:compconsistXblock}, let $R_i\in B_{X,\Phi}(G_i)$.  So $R_i=V_i\setminus A_i$ for some $A_i\subseteq V_i$ such that $A_i=\varphi(A_i)$ and $A_i$ is not an $X$-set of $G_i$.  Since $X$ is component consistent, $A_i$ is not an $X$-set of $G$.  Let $A= A_i\du (\du_{j\ne i} V_j)$.  Since $\Phi$ is component consistent, 
\[
\varphi(A)= A_i\du \lp\bigsqcup_{j\ne i} V_j\rp\mbox{ so } V\setminus \varphi(A)= (V_i\setminus A_i)\du \lp\bigsqcup_{j\ne i} V_j\setminus V_j\rp = R_i.
\]
Thus $R_i\in B_{X,\Phi}(G)$.

For \eqref{26c1}, let $R\in B_{X,\Phi}(G)$, so $R=V\setminus A$ where $A=\varphi(A)$ and $A$ is not an $X$-set.  Let $A_i=A\cap V_i$ for $i=1,\dots,k$.   Thus $R\cap V_i=(V\setminus A)\cap V_i=V_i\setminus A_i$.  Since $\Phi$ is component consistent,  $\varphi(A_i)=A_i$. Thus if $R_i=R\cap V_i\ne \emptyset$, then $A_i$ is not an $X$-set and $R_i\in B_{X,\Phi}(G_i)$. 
\epf

\begin{prop}\label{p:BXPhi-irred}
Let $X$ be a super $X$-set parameter and let $\Phi$ be an $X$-compliant closure family. Then $\B_X^{\min} \subseteq \B_{X,\Phi}$ and hence $\B_{X,\Phi}$ is an $X$-irredundant blocking family.
\end{prop}
\begin{proof}
Let $G$ be a graph with closure operator $\varphi \in \Phi$, let $R\in B_X^{\min}(G)$ and $A = V(G) \setminus R$. It suffices to show that $\varphi(A) = A$ as this implies $R\in \B_{X,\Phi}$. Since $A$ is not an $X$-set, $V(G)\setminus \varphi(A) \in B_{X,\Phi}(G)$. As $\varphi$ is a closure operator, $A\subseteq \varphi(A)$ and hence $V(G)\setminus \varphi(A) \subseteq V(G)\setminus A = R$. The minimality of $R$ implies $R =V(G)\setminus \varphi(A)$. Thus, $\varphi(A) = A$. 
 \end{proof}

\begin{cor}\label{Phi-robust}
 Let $X$ be a robust super $X$-set parameter and let $\Phi$ be an $X$-compliant inclusive component  consistent closure family.  Then $Y_{\B_{X,\Phi}}$ is robust.
\end{cor}
\bpf By Propositions   \ref{o:clos-compconst} and \ref{p:BXPhi-irred} and Lemma \ref{l:xirred-compconsist}, $Y_{\B_{X,\Phi}}$ is component consistent. Let $G$ be a connected graph of order at least two. By definition of inclusive, $y$ is in some $V(G)\setminus \varphi(A)\in B_{\Phi,x}(G)$ for every $y\in V(G)$.
Thus $\{y\}$ an $X$-irredundant set and  $Y_{\B_{X,\Phi}}$ is robust.
\epf 

 We now provide a natural description of standard forts using closure operators for standard zero forcing. 

\begin{defn}
For standard zero forcing, define $\Phi_Z = \{\varphi_{Z,G} : G \text{ is a graph}\}$, where $\varphi_{Z,G}(A)$ is the final coloring of the set $A\subseteq V(G)$ in $G$.
\end{defn}

\begin{obs}\label{o:BZ-BphiZ}
Every $\varphi_{\Z,G} \in \Phi_{\Z}$ is a closure operator and $\Phi_Z$ is a $Z$-compliant component consistent inclusive family. The blocking family already defined in Section \ref{s:irred} is the same as that defined by $\Phi_Z$, i.e., $\B_{\Z} = \B_{\Z,\Phi}$. Equivalently, $B_{\Z,\Phi}(G)$ is the set of forts of $G$. Thus $Y_{\B_{\Z,\Phi}}=\ZIR$.
\end{obs}

Like zero forcing, PSD forcing has a natural closure operator.
\begin{defn}
For PSD forcing, $\Phi_{\Zp} = \{\varphi_{\Zp,G} : G \text{ is a graph}\}$, where $\varphi_{\Zp,G}(S)$ is the final PSD coloring of the set $S$ in $G$.  Define  $B_{\Zp,\Phi}(G) = \{ V(G) \setminus \varphi_{\Zp,G}(A) :  A \text{ is not an $\Zp$-set}\}$.
\end{defn}

\begin{obs}
    Every $\varphi_{\Zp,G} \in \Phi_{\Zp}$ is a closure operator and $\Phi_{\Zp}$ is a $\Zp$-compliant  component consistent  inclusive family. Furthermore, $B_{\Zp}(G)=B_{\Zp,\Phi}(G)$, i.e., $B_{\Zp,\Phi}(G)$ is the set of PSD forts of $G$. Thus $Y_{\B_{\Zp,\Phi}}=\ZpIR$.
\end{obs}

\begin{defn}
For skew forcing, $\Phi_{\Zs} = \{\varphi_{\Zs,G} : G \text{ is a graph}\}$, where $\varphi_{\Zs,G}(S)$ is the final skew coloring of the set $S$ in $G$.  
\end{defn}

\begin{obs}
    Every $\varphi_{\Zs,G} \in \Phi_{\Zs}$ is a closure operator and $\Phi_{\Zs}$ is a $\Zs$-compliant  component consistent family.  
     Furthermore, $B_{\Zs}(G)=B_{\Zs,\Phi}(G)$, i.e., $B_{\Zs,\Phi}(G)$ is the set of skew forts of $G$. Thus $Y_{\B_{\Zs,\Phi}}=\ZsIR$.   Note that $\Phi_{\Zs}$ is not inclusive. 
\end{obs}

 It is arguable that for variants of zero forcing such as standard, PSD and skew forcing,  each of which is a propagation process, the closure operator is the  natural way to identify an appropriate blocking family.  Since the study of irredundance began  with domination, with private neighbors and thus blocking sets as neighborhoods,  we have followed that model. However, it is possible to define domination irredundant sets using closure operators.
Before we describe this, we require some { additional basic definitions and results. A family of sets is \emph{union-closed} provided the union of every pair of sets in the family is contained in the family.   It has been observed that a union of standard forts is a standard fort \cite{ZIR}, so  $B_{\Z,\Phi}(G)$  is union-closed  for every graph $G$; the same is true for skew forts \cite{skewTErecon}.  It is also true that a union of PSD forts is a PSD fort (this is slightly more complicated to see but is a consequence of the next result).

\begin{prop}\label{prop: blocking union closed}
Let $X$ be a super $X$-set parameter, let $\Phi$ be an $X$-compliant closure family. Then $B_{X,\Phi}(G)$ is a union-closed family for every $B_{X,\Phi}(G) \in \B_{X,\Phi}$.
\end{prop}
\begin{proof}
Let $G$ be a graph with closure operator $\varphi\in \Phi$, let $V=V(G)$, and let $R_1,R_2\in B_{X,\Phi}(G)$. Using the characterization in Observation \ref{o:Xcomply} and elementary set theory, 
\[\varphi(V\setminus (R_1\cup R_2)=\varphi((V\setminus R_1) \cap (V\setminus R_2))\subseteq \varphi(V\setminus R_1) \cap \varphi(V\setminus R_2)=(V\setminus R_1) \cap (V\setminus R_2)=V\setminus(R_1\cup R_2).\]
Thus $\varphi(V\setminus (R_1\cup R_2)=V\setminus(R_1\cup R_2)$ because  $A\subseteq \varphi(A)$ for every $A$.  So $R_1\cup R_2\in B_{X,\Phi}(G)$.
\end{proof}

A \emph{generator} of a union-closed family $\mathcal{F}$ is a set $A\in\mathcal{F}$ such that $A$ is not the union of two elements in $\mathcal{F}$.  Let $X$ be a super $X$-set parameter and $\Phi$ be an $X$-compliant closure family. We write $B_{X,\Phi}^*(G)$ for the set of generators of $B_{X,\Phi}(G)\in \B_{X,\Phi}$. The \emph{family of generators} for $\B_{X,\Phi}$ is defined to be $\B_{X,\Phi}^* = \{B_{X,\Phi}^*(G) : G \text{ is a graph}\}$.

\begin{obs} 
Let $X$ be a super $X$-set parameter and $\Phi$ be an $X$-compliant closure family. Then $\B_{X,\Phi}^*$ is an $X$-blocking family.
\end{obs}

\begin{obs} 
Let $X$ be a component consistent super $X$-set parameter and $\Phi$ be an $X$-compliant component consistent closure family. Then for a graph $G$, any set $R\in B_{X,\Phi}^*(G)$ satisfies $R\subseteq V(G_i)$ for some component $G_i$ of $G$.
\end{obs}

\begin{prop}\label{p:XIr-phiXIr}
Let $X$ be a super $X$-set parameter and $\Phi$ be an $X$-compliant closure family. Let $G$ be a graph and $S\subseteq V(G)$. Then $S$ is an $X$-irredundant set for $\B_{X,\Phi}$ if and only if $S$ is an $X$-irredundant set for $\B_{X,\Phi}^*$. 
Furthermore, $Y_{\B^*_{ X,\Phi}}=Y_{\B_{X,\Phi}}$.
\end{prop}
\begin{proof}
Observe that $\B_{X,\Phi}^* \subseteq \B_{X,\Phi}$. Thus, if $S$ is an $X$-irredundant set for $\B_{X,\Phi}^*$, then $S$ is an $X$-irredundant set for $\B_{X,\Phi}$. 

Suppose $S$ is an $X$-irredundant set for $\B_{X,\Phi}$. Then every $u\in S$ has a private $\B_{X,\Phi}$-set $R_u$ in $G$. If $R_u$ is a generator, then $R_u$ is a private $\B_{X,\Phi}^*$-set for $u$. So, suppose $R_u$ is not a generator. Then $R_u = R_1 \cup R_2$ for some $R_1\in B_{X,\Phi}^*(G)$ and $R_2\in B_{X,\Phi}(G)$, where $u\in R_1$. Since $R_1\subseteq R_u$, $R_1$ is a private $\B_{X,\Phi}^*$-set for $u$. Thus $S$ is an $X$-irredundant set for $\B_{X,\Phi}^*$.

The last statement is now immediate. \end{proof}

\begin{rem}
Let $X$ be a super $X$-set parameter and $\Phi$ be an $X$-compliant closure family. If $S\subseteq V(G)$ such that $S\cap R\ne\emptyset$ for all $R\in B^*_{X,\Phi}(G)$, then $S\cap R\ne\emptyset$ for all $R\in B_{X,\Phi}(G)$ from the definition of generator.  Thus $\B_{X,\Phi}^*$ is $X$-irredundant and hence $\B_X^{\min}\subseteq \B_{X,\Phi}^*$. 
\end{rem}

irredundance was defined, namely 
Recall that for domination we defined the $D$-blocking family $\B_D = \{B_D(G) : G \text{ is a graph}\}$, where $B_D(G) =\{N[v]: v\in V(G)\}$. 

\begin{defn} 
For domination, define $\Phi_D = \{\varphi_{D,G} : G \text{ is a graph}\}$, where $\varphi_{D,G}$ is given by 
\[
\varphi_{D,G}(A)=\{v:N[v]\subseteq N[A]\}  .
\]
\end{defn}

 \begin{ex}
    Let $G$ be the Bull graph shown in Figure \ref{f:bull}. 
   Then $ B_{D}(G)=\{\{1,2\},\{1,2,3,4\},\{2,3,4\},$ $\{2,3,4,5\},\{4,5\}\}$ (the set of closed neighborhoods of single vertices).  
    By computing $\varphi_{D,G}(A)$ for  $A\subseteq  V(G)$ and $|A|< 4=|V(G)|-1$, we see that $ B_{D,\Phi}(G)=\{\{1,2\},\{4,5\},\{2,3,4\}, \{1,2,3,4\},\{1,2,4,5\},\{2,3,4,5\}\}$. 
    Examination of $ B_{D,\Phi}(G)$ shows that $ B^*_{D,\Phi}(G)=\{\{1,2\},\{4,5\},\{2,3,4\}\}$.
\begin{figure}[h!]
\centering
\scalebox{1}{
\begin{tikzpicture}[scale=2,every node/.style={draw,shape=circle,outer sep=2pt,inner sep=1pt,minimum size=.2cm}]		
\node[fill=none]  (1) at (-1,1) {$1$};
\node[fill=none]  (2) at (0,1) {$2$};
\node[fill=none]  (3) at (0.5,0.14) {$3$};
\node[fill=none]  (4) at (1,1) {$4$};
\node[fill=none]  (5) at (2,1) {$5$};

\draw[thick] (1)--(2)--(3)--(4)--(5);
\draw[thick] (2)--(4);
\end{tikzpicture}}

\caption{\label{f:bull} The Bull graph}
\end{figure}
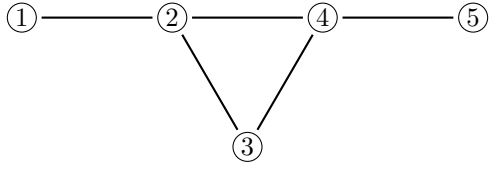\end{ex}  

The next result shows $\Phi_D$ is a  closure operator for domination and describes some elementary properties of $\Phi_D$  and its blocking sets.
 \begin{prop} 
 Let $G$ be a graph, let $V=V(G)$ and let $\varphi=\varphi_{D,G}$. 
 \ben[$(1)$]
 \item\label{dom-bas-1} Every $\varphi \in \Phi_D$ is a closure operator and $\Phi_D$ is $D$-compliant, component consistent   and inclusive.

\item\label{dom-bas-2}  For every  $A\subseteq V(G)$, $\varphi(A)= N[A]  \setminus N[V\setminus N[A]]$. 
 
\item \label{dom-bas-3} 
 $ B_D(G) \subseteq B_{D,\Phi}(G)$. 

 \item \label{dom-bas-4} 
If $W\in B_{D,\Phi}(G)$, then $
W=\cup_{w\in W, w\not\in N[V\setminus W]}N[w].$

 \item \label{dom-bas-5} 
 $B_{D,\Phi}^*(G) \subseteq B_D(G) $. 
 \een
\end{prop}

\bpf  
\eqref{dom-bas-1}:It is immediate that  $A\subseteq \varphi(A)$, $A_1\subseteq A_2$ implies $\varphi(A_1)\subseteq \varphi(A_2)$, and that  $ \Phi_{D}$ preserves isomorphism.   To see that $\varphi(\varphi(A))=\varphi(A)$, note that  $N[\varphi(A)]=N[A]$ by  definition of $\varphi(A)$.  If $y\in \varphi^2(A)$, then $N[y]\subseteq N[\varphi(A)]=N[A]$, so  $y \in \varphi(A)$.   It is immediate that $A\subseteq G_i\Rightarrow \varphi(A)\subseteq G_i$ for a component $G_i$ of $G$, and   $\Phi$ is inclusive because $\varphi(\emptyset)=\emptyset$. Finally,   $S$ is an dominating set of $G$ if and only $N[S]=V(G)$ if and only if $\varphi(S) = V(G)$.

\eqref{dom-bas-2}: For $A\subseteq V$, $v\in \varphi(A)$ if and only if $N[v]\subseteq N[A]$ if and only if $N[v]\cap (V\setminus N[A])=\emptyset$ if and only if $v\not\in N[V\setminus N[A]$ if and only if  $v\in N[A]\setminus N[V\setminus N[A]]$. 

\eqref{dom-bas-3}: Fix an element of $B_D(G)$, which is $N[u]$ for some $u\in V$. Let $A=V\setminus N[u]$. We show that $\varphi(A)\subseteq A$, which implies $\varphi(A)= A$  and thus  $ N[u]\in  B_{D,\Phi}(G)$ by Observation \ref{o:Xcomply}.   
Suppose $y\not\in A$.  Then $y\in N[u]$. Since $u\not\in N[A]$, $N[y]\not\subseteq N[A]$.  Thus $y\not\in \varphi(A)$, and therefore $\varphi(A)\subseteq A$.

\eqref{dom-bas-4}:  Let $W\in B_{D,\Phi}(G)$, so $W=V\setminus A$ for some $A=\varphi(A)$.
 Since $\varphi(A)= N[A]  \setminus N[V\setminus N[A]]$, $W=V\setminus \varphi(A) =N[V\setminus N[A]]$.  
 Thus 
$W=\cup_{w\in W, w\not\in N[V\setminus W]}N[w].$ 

\eqref{dom-bas-5}: This result is immediate from \ref{dom-bas-4}.
\epf

 In a similar fashion to domination, blocking sets for vertex covering can be obtained from a $VC$-compliant closure operator.

\begin{defn}
For vertex covering, define $\Phi_{VC} = \{\varphi_{VC,G} : G \text{ is a graph}\}$, where $\varphi_{VC,G}$ is given by 
\[
\varphi_{VC,G}(A)= A \cup\{v:N(v)\subseteq A\}.
\]
\end{defn}
 Observe that if $G$ has isolated vertices $v_1,\dots,v_r$, then $\varphi_{VC}(\emptyset)=\{v_1,\dots,v_r\}$.
 
 \begin{ex}
    Let $G$ be the Paw graph shown in Figure \ref{f:paw}. 
   Then $ B_{VC}(G)=\{\{1,4\},\{2,3\},\{2,4\},\{3,4\}\}$ is the set of edges of $G$.  
    By computing $\varphi_{VC,G}(A)$ for $A\subseteq  V(G)$ and $|A|< 3=|V(G)|-1$, we see that $ B_{VC,\Phi}(G)=\{\{1,4\},\{2,3\},\{2,4\},\{3,4\}\},\{1,2,4\}, \{1,3,4\}, \{2,3,4\}$. 
    Examination of $ B_{VC,\Phi}(G)$ shows that $ B^*_{VC,\Phi}(G)= B_{VC}(G)$ (this is true in general, as seen in the next result).
\end{ex}  
  
  \begin{prop} 
 Let $G$ be a graph, let $V=V(G)$ and let $\varphi=\varphi_{VC,G}$.
 \ben[$(1)$]
 \item\label{vc-bas-1} Every $\varphi \in \Phi_{VC}$ is a closure operator and $\Phi_{VC}$ is $VC$-compliant, component consistent,  and inclusive.

\item \label{vc-bas-3} 
 $ B_{VC}(G) \subseteq B_{VC,\Phi}(G)$. 

 \item \label{vc-bas-4} If $W\in B_{VC,\Phi}(G)$, then  $W$ is a union of edges.

 \item \label{vc-bas-5} 
 $ B_{VC,\Phi}^*(G)  = B_{VC}(G) $. 
 \een
\end{prop}

\bpf  
\eqref{vc-bas-1}: Let $A, A_1, A_2\subseteq V$. It is immediate that  $A\subseteq \varphi(A)$, $A_1\subseteq A_2$ implies $\varphi(A_1)\subseteq \varphi(A_2)$, and that  $ \Phi_{VC}$ preserves isomorphism. If $v\not\in \varphi(A)$, then  $v\not\in A$ and $v$ is adjacent to $w\not\in A$, which implies $w\not\in\varphi(A)$ and thus $v,w\not\in\varphi(\varphi(A))$.  Hence $\varphi(\varphi(A))=\varphi(A)$.  
 It is immediate that $A\subseteq G_i\Rightarrow \varphi(A)\subseteq G_i$ for a component $G_i$ of $G$, and   $\Phi_{VC}$ is inclusive because if $G$ is connected and of order at least two, then $\varphi(\emptyset)=\emptyset$.
To see that $\Phi_{VC}$ is $VC$-compliant:
 $S$ is a $VC$-set if and only if every edge has an endpoint in $S$ if and only if for all $v\in V$, $v\in S$ or $N(v)\subseteq S$ if and only if $\varphi(S)=V$. 

\eqref{vc-bas-3}:  Recall that $B_{VC}(G)$ is the set $E$ of edges of $G$.  Let $\{v,w\}\in E$.  Then $v,w\not\in\varphi(V\setminus \{v,w\})$, so $\{v,w\}\in B_{VC,\Phi}$.

\eqref{vc-bas-4}:  Let $W\in B_{VC,\Phi}(G)$, so $W=V\setminus A$ for some $A=\varphi(A)$.  Let $w\in W$.  Since $w\not\in\varphi(A)$, $w\not\in A$ and $N(w)\not\subseteq A$.  Thus, there is some $u\in N(w)$ such that $u\not\in A$.  Thus, $\{w,u\}\subseteq W$.

\eqref{vc-bas-5}: This result is immediate from \ref{vc-bas-4}.
\epf

There is a well-known conjecture, due to Peter Frankl, about union-closed families: For a finite union-closed family of sets (other than the family containing just the empty set) there exists an element that is contained in at least half of the sets in the family. A proof of the union-closed conjecture has been claimed (see \cite{unionclosedconj}) but has not yet been published.

\begin{rem} If the union-closed sets conjecture is true, then for every graph $G$, there exists a vertex $v\in V(G)$ that is contained in at least half the (standard or PSD) forts of $G$.  \end{rem}

This does not apply to skew forts, because $B_{\Zs}(G)=\emptyset$ when $\Zs(G)=0$.  Furthermore, the domination blocking family $B_D(G)=\{N(v):v\in V(G)\}$ and  the vertex cover blocking family $B_{VC}(G)=E(G)$ are not union closed and the path $P_n,n\ge 7$ is an example where no vertex is in at least have of the blocking sets.

\section{Concluding remarks}\label{s:conclude}
 In this section we summarize the values of lower XIr number, $X$ number, upper $X$ number, and upper XIr number for the graph parameters standard zero forcing, PSD forcing, skew forcing, and vertex covering, for various graph families (with each graph family in a separate table).  
 Open questions for future work include Questions \ref{q:delta-Z+} and \ref{q:Zs-forrest}.

The complete graphs  and empty graphs   do not need tables for their summaries since they are almost all the same: For the complete graph, $\zir(K_n)=\Z(K_n)=\zbar(K_n)=\ZIR(K_n)=\zpir(K_n)=\Zp(K_n)=\zpbar(K_n)=\ZpIR(K_n)=\vc(K_n)=\vcbar(K_n)=\VCIR(K_n)=n-1$ and $\zsir(K_n)=\Zs(K_n)=\zsbar(K_n)=\ZsIR(K_n)=n-2$. 
For the empty graph, $\zir(\ol{K_n})=\Z(\ol{K_n})=\zbar(\ol{K_n})=\ZIR(\ol{K_n})=\zpir(\ol{K_n})=\Zp(\ol{K_n})=\zpbar(\ol{K_n})=\ZpIR(\ol{K_n})=\zsir(\ol{K_n})=\Zs(\ol{K_n})=\zsbar(\ol{K_n})=\ZsIR(\ol{K_n})=n$ and $\vc(\ol{K_n})=\vcbar(\ol{K_n})=\VCIR(\ol{K_n})=0$.
Table \ref{tableKpq} presents complete bipartite graphs. Path and cycles are presented in  Tables  \ref{tablePn} and \ref{tableCn}. For each of these graphs,  the value of $X(G)$ for these four parameters has been known for years; 
these values and the values of $\xbar(G)$   appear in \cite{XTAR-YTAR}.

\begin{table}[h!]
\renewcommand{\arraystretch}{1.2}
\begin{center} 
\noindent { \begin{tabular}{| l  || c | c| c| c| c| c|}
\hline
  $K_{p,q}$ &  $\xir(G)$ & $X(G)$ & $\xbar(G)$ & $\XIR(G)$ & result \#\\
\hline\hline
$\Z$   &  $p$ & $q+p-2$  & $q+p-2$ & $q+p-2$ & \cite{ZIR} \\
\hline
$\Zp, p\ge 2$  &  $p$ & $p$  & $q$ & $\max(q,q+p-4)$ & \ref{ex-star-PSD} \\
$\null\ \ \ \ \ p=1$  &  1 & 1  & 1 & 1 & \ref{ex:PSDtree} \\\hline
$\Zs$   &  $q+p-2$ & $q+p-2$  & $q+p-2$ & $q+p-2$ & \ref{ex-star-skew} \\
\hline
$\vc,\ p\ge 2$   &  $p$ & $p$  & $q$ & $q+p-2$ & \ref{e:VCIRbigger} \\
$\null\ \ \ \ p=1$  &  1 & 1  & q & q &  \\\hline
\end{tabular}}
\caption{Irredundance-related parameter values for the complete bipartite graph $K_{p,q}$ \vspace{-15pt}}\label{tableKpq}
\end{center}
\end{table}
\begin{table}[h!]
\renewcommand{\arraystretch}{1.4}
\begin{center} 
\noindent { \begin{tabular}{| l  || c | c| c| c| c| c|}
\hline
  $P_n$ &  $\xir(G)$ & $X(G)$ & $\xbar(G)$ & $\XIR(G)$ & result \#\\
\hline\hline
$\Z,  n\ge 5$   &  $1$ & $1$  & $2$ & $\lf\frac{n-1}2\rf$ & \cite{ZIR} \\
\hline
$\Zp$  &  $1$ & $1$  & $1$ & $1$ &  \ref{ex:PSDtree} \\
\hline
$\Zs$, $n$ even   &  $0$ & $0$  & $0$ & $0$ & \ref{ex:path-skew} \\
$\null \ \ \ \ \ \ n$ odd &  $1$ & $1$  & $1$ & $1$ &  \\
\hline
$\vc$   &  $\lf\frac n 2\rf$ & $\lf\frac n 2\rf$  & $\lc\frac {2n-1} 3\rc$ & $\lc\frac {2n-1} 3\rc$ & \ref{ex:VC-path-cyc}\\
\hline
\end{tabular}}
\caption{Irredundance-related parameter values for the path $P_n$ \vspace{-15pt}
}\label{tablePn}
\end{center}
\end{table}
\begin{table}[h!]
\renewcommand{\arraystretch}{1.4}
\begin{center} 
\noindent { \begin{tabular}{| l  || c | c| c| c| c| c|}
\hline
  $C_n,   n\ge 4$ &  $\xir(G)$ & $X(G)$ & $\xbar(G)$ & $\XIR(G)$ & result \#\\
\hline\hline
$\Z$   &  $2$ & $2$  & $2$ & $\lf\frac{n}2\rf$ & \cite{ZIR} \\
\hline
$\Zp$  &  $2$ & $2$  & $2$ & $2$ &  \ref{Zp-cycle} \\
\hline
$\Zs$, $n$ even   &  $2$ & $2$  & $2$ & $2$ & \ref{ex-Cn-skew} \\
$\null \ \ \ \ \ \ n$ odd &  $1$ & $1$  & $1$ & $1$ &  \\
\hline
$\vc$   &  $\lc\frac n 2\rc$ & $\lc\frac n 2\rc$  & $\lc\frac {2n-2} 3\rc$ & $\lc\frac {2n-2} 3\rc$ & \ref{ex:VC-path-cyc} \\
\hline
\end{tabular}}
\caption{Irredundance-related parameter values for the cycle $C_n$ 
}\label{tableCn}
\end{center}
\end{table}



\section*{Statements and Declarations}
The research of B.~Curtis  was partially supported by National Science Foundation grant  1839918.  The research of M.~Flagg  is partially supported by National Science Foundation grant  2331634.  This research began at the American Institute of Mathematics (AIM) and the authors thank AIM and  NSF for their support.

Data used to find patterns and examples was generated by SAGE code. Available as code for download as Jupyter notebook and  PDF  at  \url{https://aimath.org/~hogben/Uni-Irred-ZpZsZ.ipynb}, \url{https://aimath.org/~hogben/Uni-Irred-ZpZsZ.pdf}. 

This paper was approved for public release: AFRL-2025-4326. The views expressed in this article are those of the authors and do not necessarily reflect the official policy or position of the Air Force, the Department of Defense, or the U.S. Government.


\begin{thebibliography}{1}


\bibitem{AIM08} AIM Minimum Rank -- Special Graphs Work Group (F. Barioli, W. Barrett, S. Butler,   S. M. Cioaba, D. Cvetkovi\'c,  S. M. Fallat, C. Godsil,  W. Haemers,  L. Hogben,  R. Mikkelson,  S. Narayan,  O. Pryporova,   I. Sciriha,  W. So,   D. Stevanovi\'c,  H. van der Holst,K. Vander Meulen, and A. Wangsness).  Zero forcing sets and the minimum rank  of graphs.   {\em Lin. Alg. Appl.} 428 (2008), 1628--1648.

\bibitem{smallparam} F. Barioli, W. Barrett, S.M. Fallat, H.T. Hall, L. Hogben, B. Shader, P. van den Driessche, and H. van der Holst.  Zero forcing parameters and minimum rank problems. {\it{Linear Algebra Appl.}} {433} (2010), 401-411.


\bibitem{BBCF20domchain}
C.~Bazgan, L.~Brankovic, K.~Casel, and H.~Fernau.  Domination chain: characterisation, classical complexity, parameterised complexity and approximability.
\emph{Discrete Appl. Math.} 280 (2020), 23--42.


\bibitem{PDrecon} B. Bjorkman, C. Bozeman, D. Ferrero, M. Flagg, C. Grood, L. Hogben, B. Jacob, and C. Reinhart.  Power domination reconfiguration.  To appear in \emph{La Matematica}.  Available at \url{https://arxiv.org/abs/2201.01798}.

\bibitem{skewTErecon} N.~Bong, M.~Flagg, M.~Hunnell, J.~Hutchins, R.~Moruzzi, H.~Schuerger, and B.~Small. Reconfiguration of Minimum PSD Forcing Sets and Minimum Skew Forcing Sets. Under review.  Available at \url{https://arxiv.org/abs/2501.03642}.

\bibitem{XTARiso} N.H. Bong, J. Carlson, B. Curtis, R. Haas, and L. Hogben. Isomorphisms  and  properties of TAR reconfiguration graphs for zero forcing and other $X$-set parameters.  \emph{Graphs Combin.} {\bf 39} (2023), Paper No. 86, 23 pp. 

\bibitem{BC22}  B.~Brimkov and J.~Carlson. Minimal Zero Forcing Sets. \emph{
Australas.~J.~Combin.} 90 (2024), 363–377.

\bibitem{BFH19} B.~Brimkov, C.C.~Fast, and I.V.~Hicks. Computational approaches for zero forcing and related problems. \emph{European J. Oper. Res.} 273 (2019),  889--903.


\bibitem{CKMU19hyper} 
A.~Conte, M.M.~Kanté, A.~Marino, amd T.~Uno.
Maximal irredundant set enumeration in bounded-degeneracy and bounded-degree hypergraphs.
In \emph{Combinatorial
Algorithms}, C.J.~Colbourn,
R.~Grossi, and 
N.~Pisanti, Eds., \emph{Lecture Notes in Comput. Sci.}, 11638
Springer, Cham, 2019, pp. 148–159.

\bibitem{CHM78}  E.J.~Cockayne, S.T.~Hedetniemi, and D.J.~Miller.
Properties of hereditary hypergraphs and middle graphs.
\emph{Canad. Math. Bull.}, 21 (1978), 461--468.

\bibitem{CM93} E.J.~Cockayne and C.M.~Mynhardt. The sequence of upper and lower domination, independence
and irredundance numbers of a graph. \emph{Discrete Math.} 122 (1993), 89--102.

\bibitem{XTAR-YTAR} B.A.~Curtis, M.K.~Flagg, and L.~Hogben. TAR reconfiguration for vertex set parameters. Under review.  Available at \url{https://arxiv.org/abs/2406.05509}.

\bibitem{sage:irred} B.~Curtis and L.~Hogben. {\em Sage} code for irredundance for PSD forcing and skew forcing. 
Available  as code for download as Jupyter notebook and  PDF  at  \url{https://aimath.org/~hogben/Uni-Irred-ZpZsZ.ipynb}, \url{https://aimath.org/~hogben/Uni-Irred-ZpZsZ.pdf}. 




\bibitem{ZIR} B.~Curtis, L.~Hogben, and A.~Roux. Zero forcing irredundant sets. Under review.  Available at \url{https://arxiv.org/abs/2403.03921}.


\bibitem{unionclosedconj}
R.~Demontis. The union-closed set conjecture is true.  Available at \url{https://arxiv.org/abs/2405.03731}.



\bibitem{FH18}
{C.C.~Fast and I.V.~Hicks}, {Effects of vertex degrees on the
  zero-forcing number and propagation time of a graph}, {\em Discrete Appl.
  Math.}, 250 (2018), pp.~215--226.


\bibitem{FHHHK02} 
O.~Favaron, T.W.~Haynes, S.T.~Hedetniemi, M.A.~Henning, and D.J.~Knisley. Total irredundance in graphs. 
\emph{Discrete Math.} 256 (2002),  115--127.

\bibitem{HS14} R.~Haas and K.~Seyffarth. The $k$-Dominating graph. \emph{Graphs and Combin.}, 30 (2014) 609--617. 



\bibitem{HHS98} T.W.~Haynes, S.T.~Hedetniemi, and P.J.~Slater.
\emph{Fundamentals of domination in graphs}. 
Marcel Dekker, Inc., New York, 1998.


 \bibitem{HLS22book} L.~Hogben, J.C.-H.~Lin, and B.L.~Shader. {\em Inverse Problems and Zero Forcing for Graphs}. {Mathematical Surveys and Monographs} {270}, American Mathematical Society, Providence, RI,  2022. 

 
\bibitem{IMA10} IMA-ISU research group on minimum rank (M. Allison,   E. Bodine, L. M.  DeAlba, J. Debnath, L. DeLoss, C. Garnett, J. Grout, L. Hogben, B. Im, H. Kim, R. Nair, O. Pryporova, K. Savage, B.  Shader, and A. Wangsness Wehe).  Minimum rank of skew-symmetric matrices described by a graph.   {\em Lin. Alg. Appl.} 432 (2010), 2457--2472.


\bibitem{King15} N.F.~Kingsley.  Skew propagation time. Dissertation (Ph.D.), Iowa State  University, 2015.

\bibitem{MR21} C.M.~Mynhardt and A.~Roux. Irredundance. In \emph{Structures of domination in graphs} (Dev. Math., 66), T.W.~Haynes, S.T.~Hedetniemi, and M.A.~Henning, Editors, Springer, Cham, 2021, pp.  135--181.


\bibitem{SMH19_PSDforts} L.A.~Smith, D.J.~Mikesell, and I.V.~Hicks.  An integer program for positive semidefinite zero forcing in graphs.
\emph{Networks}, 76 (2020), 366--380.


\bibitem{T93boolean} S.~Todor\v cevi\' c. Irredundant sets in Boolean algebras. \emph{Trans. Amer. Math. Soc.} 339 (1993), 35--44.
\end{thebibliography}
\end{document}